\documentclass[opre,nonblindrev]{informs3_revised2}

\OneAndAHalfSpacedXI 

\usepackage{endnotes}
\let\footnote=\endnote
\let\enotesize=\normalsize
\def\notesname{Endnotes}%
\def\makeenmark{$^{\theenmark}$}
\def\enoteformat{\rightskip0pt\leftskip0pt\parindent=1.75em
  \leavevmode\llap{\theenmark.\enskip}}

\usepackage{graphicx}
\usepackage{subcaption}
\usepackage{mwe}
\usepackage{amssymb}
\usepackage{multirow}

\def\makeenmark{$^{\theenmark}$}
\def\enoteformat{\rightskip0pt\leftskip0pt\parindent=1.75em
  \leavevmode\llap{\theenmark.\enskip}}

\usepackage{natbib}
 \bibpunct[, ]{(}{)}{,}{a}{}{,}%
 \def\bibfont{\small}%
 \def\bibsep{\smallskipamount}%
 \def\bibhang{24pt}%
 \def\newblock{\ }%
\TheoremsNumberedThrough     
\ECRepeatTheorems

\EquationsNumberedThrough    

\begin{document}


\RUNAUTHOR{Sasanuma, Hampshire, and Scheller-Wolf}

\RUNTITLE{Markov Chain Decomposition Based On Total Expectation Theorem}

\TITLE{Markov Chain Decomposition Based On Total Expectation Theorem}

\ARTICLEAUTHORS{%
\AUTHOR{Katsunobu Sasanuma}
\AFF{College of Business, Stony Brook University, Stony Brook, NY 11794, \EMAIL{katsunobu.sasanuma@stonybrook.edu}} 
\AUTHOR{Robert Hampshire}
\AFF{Gerald R. Ford School of Public Policy, University of Michigan, Ann Arbor, MI 48109, \EMAIL{hamp@umich.edu}}
\AUTHOR{Alan Scheller-Wolf}
\AFF{Tepper School of Business, Carnegie Mellon University, Pittsburgh, PA 15213, \EMAIL{awolf@andrew.cmu.edu}}
} 
\ABSTRACT{\textcolor{black}{A divide-and-conquer approach to analyzing Markov chains (MCs) is not utilized as widely as it could be, despite its potential benefits. One primary reason for this is the fact that most MC decomposition approaches involve a complex and inflexible methodology: decomposed subchains must be disjoint, transition rates of these decomposed subchains must be altered in a way tailored to the particular MC model, and the procedure to aggregate suchains needs to incorporate a nonlinear normalization constraint, complicating the analytical expression of performance measures. In contrast, we propose a versatile yet simple decomposition method for continuous time MCs based on the total expectation theorem. Leveraging the properties of this theorem, our method has great flexibility in the choice of subchains, and the procedure to obtain expected values of interest is simply a linear summation of subchains' properties, which is not affected by the normalization constraint. We prove that to maintain the correct distribution of decomposed subchains one may use our novel \emph{termination} scheme, a modification of transition rates, that ensures \emph{partial flow conservation} at boundary states. This termination scheme is applicable to MCs with any structure, since the scheme depends only on the boundary-state distribution, not on the structure of the MCs. To demonstrate the generality and capability of our method, we analytically solve various models, such as a congestion-based staffing queue and a Markov-modulated M\textsubscript{t}/M\textsubscript{t}/1 queue. As not all systems admit an analytical solution, we complement this analysis with numerical studies of MCs with various sizes using the algorithm based on our method.}
}%


\maketitle

%


\section{Introduction}
\textcolor{black}{A divide-and-conquer approach is commonly used to find a quantity of interest in a complex system. This approach is generally desirable because a ``complex" system is often composed of multiple simpler subsystems, each of which may be easy to analyze in isolation. In fact, a variety of divide-and-conquer methods are used in many fields, including the analysis of Markov chains (MCs). \cite{stewart1994introduction} reviews major MC decomposition methods that have been developed and utilized. These methods allow for an exact (or approximate) analysis of various performance measures, such as expected values of interest, which we consider in this paper. These methods are typically composed of two main steps: (1) a decomposition procedure, and (2) an aggregation procedure.}

\textcolor{black}{The first step is to decompose a full MC into a set of subchains, which is almost always a disjoint set because most methods require the decomposition of the transition matrix corresponding to an MC system. These subchains must maintain their stationary distributions after decomposition (i.e., proportional to the original full MC). This conservation of distribution is essential in MC decomposition methods; various techniques have been developed to satisfy this requirement. For example, the NCD (\emph{nearly completely decomposable}) method decomposes a full MC into a set of almost independent subchains to approximate solutions of MC models \citep{simon1961aggregation}---but subchains comprising the full MC are not always nearly independent from each other. The \emph{censoring} method, another very popular decomposition technique, creates a modified (augmented) censored transition matrix corresponding to a subchain to analyze the properties of MC systems---but in general, the task of computing the correct censored transition matrix (that conserves the distribution of the censored MC) is non-trivial \citep{freedman1983approximating, zhao1996censored}.}

\textcolor{black}{The second step is to aggregate the solutions of all decomposed subchains. The standard aggregation procedure finds probability weights for subchains by imposing a normalization condition, retrieves the distribution of the full MC using the weights, and then uses this distribution to derive the expected values of interest, e.g., average number of jobs, average waiting time, z-transforms. This standard procedure is straightforward, but the analytical representation of performance measures becomes a nonlinear function of properties of subchains, obscuring the relationship between subchains and the original full MC.}

\textcolor{black}{For many, the benefit obtained by taking a divide-and-conquer approach does not justify the cost of extra work involved in decomposition, aggregation, and normalization, inhibiting widespread usage of MC decomposition methods in general. In fact, few take a decomposition approach to solve an MC if it can be solved as is. To remedy these shortcomings, we provide a new MC decomposition method based on the total expectation theorem (Theorem \ref{total expectation}). Leveraging the total expectation theorem (the law of iterated expectations), we can represent a performance measure (an expected value of interest) as a linear combination of the performance measures of decomposed subchains, where a choice of subchains is flexible: we are allowed to use overlapping or nested subchains if such a choice is convenient for analyzing the MC. This linear representation also captures the normalization condition as a special case, greatly simplifying the application of our methodology. The specific procedure to impose the normalization condition is similar to the one used for z-transforms (or moment-generating functions): We introduce a single unknown parameter representing the reference state (or set) probability, which is easily found at the end of the analysis by setting the argument of all expectations equal to 1.}

\textcolor{black}{Our method also provides a simple way to maintain the distributions of decomposed subchains using the key idea that the steady-state distribution of a subchain remains unchanged as long as average flows are maintained. A similar idea has been utilized to solve general queues (not necessarily Markovian); see, for example, the \emph{Level Crossing Theory} (LCT) \citep{brill1977level}, the \emph{rate balance principle} (RBP) \citep{oz2017rate}, and the \emph{Queueing and Markov Chain Decomposition} (QMCD) method \citep{abouee2016state,baron2018state}. Using this key idea, we show that the steady-state distribution of a subchain is conserved if and only if the net average inflow (outflow) is conserved at all boundary states of the subchain (our \emph{partial flow conservation} condition; Lemma~\ref{partial flow}). This condition has a large degree of freedom and thus gives us great flexibility to choose a convenient termination scheme (Corollary~\ref{special termination}), which is explicitly represented by the boundary-state distribution. Thanks to the Markov property, our scheme always works regardless of the structure of MCs, just like the procedure used for the hidden Markov model (HMM).}

\subsection{Literature Review}
Over the past several decades, various MC decomposition methods have been developed. For the simplest example, if a detailed (or partial) flow balance condition holds for an MC, this MC is called reversible (or quasi-reversible), and any truncated subchains maintain the steady-state probability distribution of the full MC (up to a normalization constant). Hence, properties of these truncated subchains can be analyzed independently \citep{whittle1986systems,kellyreversibility}. Once properties of subchains are solved independently, the results are combined using a normalization condition to obtain the properties of the full MC. If an MC is ``lumpable'' \citep{kemeny1960finite}, then the MC can be partitioned into multiple subchains that maintain the steady-state probability distribution of the full MC. The method of lumping has been extended by many researchers (see, for example, \citealt{stewart1994introduction}). In a special case where a subchain has a single input state from other subchains, an MC is called a Single-Input Superstate Decomposable Markov Chain (SISDMC), which makes the lumping procedure simpler \citep{feinberg1987method}. This lumping procedure can be repeatedly applied to simplify the analysis of MCs \citep{katehakis2012successive}. However, to cope with general MCs, we need different methods.

\textcolor{black}{To decompose and analyze MCs with more general structures, the \emph{censoring} (or \emph{watching}) technique is well-known. This method is based on the \emph{censored} transition matrix, obtained by \emph{watching} the sample paths of the MC on a portion of its state space; the transitions in the complementary set are omitted. This technique has been utilized to solve advanced stochastic processes \cite[see, for example,][]{freedman1983approximating, zhao1996censored, zhao2000censoring}. The method has been further extended to analyze not only steady-state performances, but also transient properties (e.g., transient probabilities, first passage times) of more general continuous time MCs utilizing UL- and LU-types of \emph{RG}-factorizations \citep{li2010constructive}.}

\textcolor{black}{When solving general, both Markovian and non-Markovian, queues, \emph{Queueing and Markov Chain Decomposition} (QMCD) is a powerful MC decomposition technique \citep[e.g., see][]{abouee2016state,baron2018state,abouee2012strategies, wang2015m, wang2019tandem}. This method utilizes four steps \cite[see \S4 of][]{baron2018state}: 1) Decomposition of the system into subsystems 2) Tying the subsystems together, 3) Solving each subsystem, and 4) Normalizing the solution. The second step is fundamental in QMCD; its objective is to maintain the correct distribution of subsystems so that we can solve each subsystem individually; this step utilizes, for example, the supplementary variable method, the probabilistic approach, or busy/idle period analysis to maintain the properties of subsystems. Like the censoring technique, QMCD requires a careful analysis of the structure of MCs before solving each subsystem. Our MC decomposition method is similar to QMCD, except that in our method the second step utilizes termination and the fourth step utilizes the total expectation theorem, which in many cases may yield a more straightforward analysis.}

Another MC decomposition approach is the approximation method specifically developed for a \emph{nearly completely decomposable} (NCD) MC \citep{simon1961aggregation}. Under the NCD condition, a large MC is clustered into a small number of subchains, each of which is relatively independent from other subchains. By introducing a coupling matrix representing the transitions among these subchains, an approximate solution of the full MC is efficiently calculated. Lastly, if we need to analyze an MC with a general structure using a divide-and-conquer approach, recursive algorithms are often effective. One of the most popular such methods is the iterative aggregation/disaggregation (IAD) algorithm by \cite{takahashi1975lumping}, which is suitable for solving large MCs numerically. This method repeatedly decomposes a full MC into partitions (disjoint subchains) and aggregates the solutions of subchains until the equilibrium is reached. \textcolor{black}{Both the NCD approximation and the IAD method are popular decomposition methods and many variants of these algorithms have been proposed (for a review, please refer to \S 6 of \citealt{stewart1994introduction}).}

\textcolor{black}{Among the many MC decomposition methods that have been developed, our decomposition method is most related to the theory of \emph{partial (flow) balance} \citep{whittle1986systems,kellyreversibility}. Specifically, our method has a root in Theorem 9.5 of \cite{kellyreversibility}, which proves that the partial (flow) \emph{balance} condition is equivalent to conservation of the steady-state distribution between the original full MC and its truncated subchains, up to a normalization constant. We show that it is possible to decompose any MC while preserving the distribution if any partial flow \emph{imbalance} (i.e., net average flow) is conserved at all ``boundary'' states identified in the decomposition process. The required conservation of net inflow or outflow (\emph{partial flow conservation}) is made possible by adding new transitions (\emph{termination}) among the boundary states of subchains. We prove the equivalence of our partial flow conservation and the conservation of the steady-state distribution of the terminated subchain, up to a normalization constant (Lemma \ref{partial flow}). Using this result, we propose a new termination scheme (Corollary \ref{special termination}; a sufficient condition for Lemma \ref{partial flow}), with which a decomposed subchain satisfies partial flow conservation and thus, maintains the distribution of the original MC. In the special case when partial flow is balanced at all boundary states, Lemma \ref{partial flow} is reduced to Theorem 9.5 of \cite{kellyreversibility}.}

\textcolor{black}{Our decomposition method is functionally similar to the censoring technique; both methods try to maintain the correct steady-state distribution of the decomposed subchain\textemdash either by terminating the subchain in our method or by augmenting the transition matrix (the procedure to make the matrix stochastic) in the censoring method. This difference stems from the roots of these two methods: The censoring method is based on the transition matrix of the censored MC, whose sample paths are correctly maintained (after transitions at the complementary set are deleted); thus, transient properties can be well analyzed. Furthermore, the censoring procedure does not require information about the boundary-state subsets. In contrast, our method is based on partial flow conservation, which guarantees conservation of the stationary distribution of the original MC. Thus, our method can deal with steady-state properties, but not transient properties (since return probabilities to the decomposed subchain are not maintained using our termination scheme).}

\textcolor{black}{It is known that finding the correct censored transition matrix is not easy in general \citep{freedman1983approximating}. Our method overcomes this issue by utilizing the flexibility to choose any termination scheme that satisfies partial flow conservation; we choose the simplest termination scheme (Corollary \ref{special termination}) to analytically or numerically obtain the stationary probabilities of subchains on the restricted state space. Our termination scheme is easily constructed for MCs with \emph{any} structure because transitions made within the complementary set become independent given the boundary-state distribution of the complementary set; the property we use here is the Markov property, which contributes to the simplicity of our termination scheme and our decomposition method. In other words, using our termination scheme, we can always set up equations for steady-state distributions of subchains, which reveal how subchains impact each other; see, for example, \S\ref{sec:stackedqueue}, where a single model, two stacked queues, can explain how all subchains in the M\textsubscript{t}/M\textsubscript{t}/1 queue depend on each other via termination. Furthermore, our decomposition approach can effortlessly reveal hidden relationships among the full MC and its subchains (possibly overlapping with each other) thanks to the linear property of the total expectation theorem, which is not affected by a normalization condition. For example, in \S\ref{sec: CBS}, we observe how subchains A and B contribute to performance measures of the CBS model; similarly, in \cite{sasanuma2021approximate}, the authors obtain insights into the impact of each stage of a reneging queue on performance measures following our MC decomposition approach. To summarize, our method provides analytical insights into the relationships among subchains, as well as their contributions to performance measures (expectations) of interest. Alternatively, if we do not care for such relationships, we can simply decompose a system into arbitrary disjoint sets and implement our MC decomposition algorithm described in \S\ref{sec:algorithm} to obtain solutions numerically.}

\textcolor{black}{The rest of the paper is organized as follows: We first illustrate the benefits of our decomposition method using a simple example in \S\ref{sec:four-state}. We then prove the total expectation theorem in MC settings under the \emph{conservation of steady-state distribution} condition in \S\ref{sec:tet}. We present the necessary and sufficient condition to conserve the distribution of the decomposed subchain, our \emph{partial flow conservation} condition, which is satisfied by \emph{termination}, a modification of transitions in \S\ref{sec:term}. We then show three applications of our method in \S\ref{sec:appl}, providing more applications in the online appendix. \S\ref{sec:concl} concludes the paper. All proofs are presented in Appendix \ref{sec: EC-proofs} (online appendix).}
\subsection{\textcolor{black}{Example: Four-State CTMC}}
\label{sec:four-state}
\begin{figure}[h]
\FIGURE
{\includegraphics*[scale=0.33]{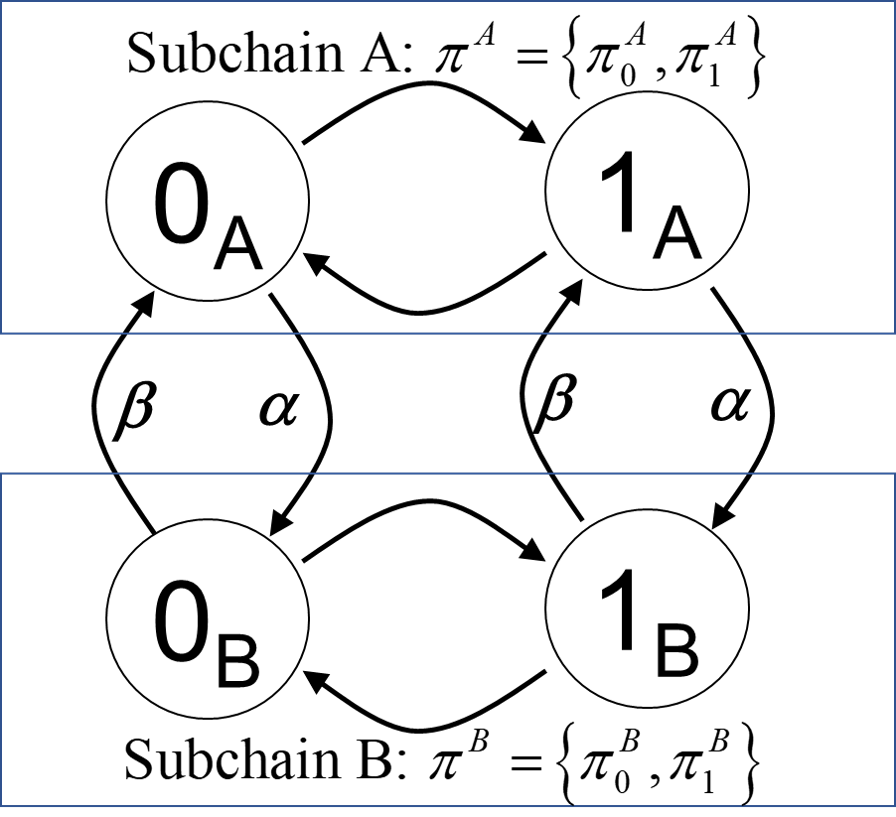}}
{Four-State CTMC.\label{fig:fourstateMC3}}
{}
\end{figure}

\textcolor{black}{Consider a simple MC with four recurrent states, which correspond to two upper states $A=\{0_A, 1_A\}$ and two lower states $B=\{0_B, 1_B\}$ (see Figure \ref{fig:fourstateMC3}). States in $A$ and states in $B$ are connected via transitions with rate $\alpha$ (rate $\beta$) from each state in $A$ ($B$) to the corresponding state in $B$ ($A$, respectively). (Note: An extension of this model is discussed in Appendix \ref{sec: extension}.) Let $\Omega=A \cup B$. Let $X$ be the random variable representing a state in $\Omega$. Define the distributions of $A$, $B$, and $\Omega$ as $\pi^A=(\pi^A_0,\pi^A_1)$, $\pi^B=\{\pi^B_0,\pi^B_1\}$, and $\pi=(\pi_{0A},\pi_{1A},\pi_{0B},\pi_{1B}\}$, respectively. Suppose that our goal is to find a performance measure $E[f(X)]\doteq\sum_{i \in \Omega} f(i) \pi_i$, where $f(X)$ is any function of $X$. For this small model, most people would solve the problem without resorting to any decomposition methods because a direct approach would be sufficient if we are only interested in a solution; however, our decomposition method provides additional information about the model, such as how subchains are related to each other and how each subchain contributes to the performance measure. We explain our method in two steps, decomposition and aggregation.}

\begin{figure}[h]
\FIGURE
{\includegraphics*[scale=0.33]{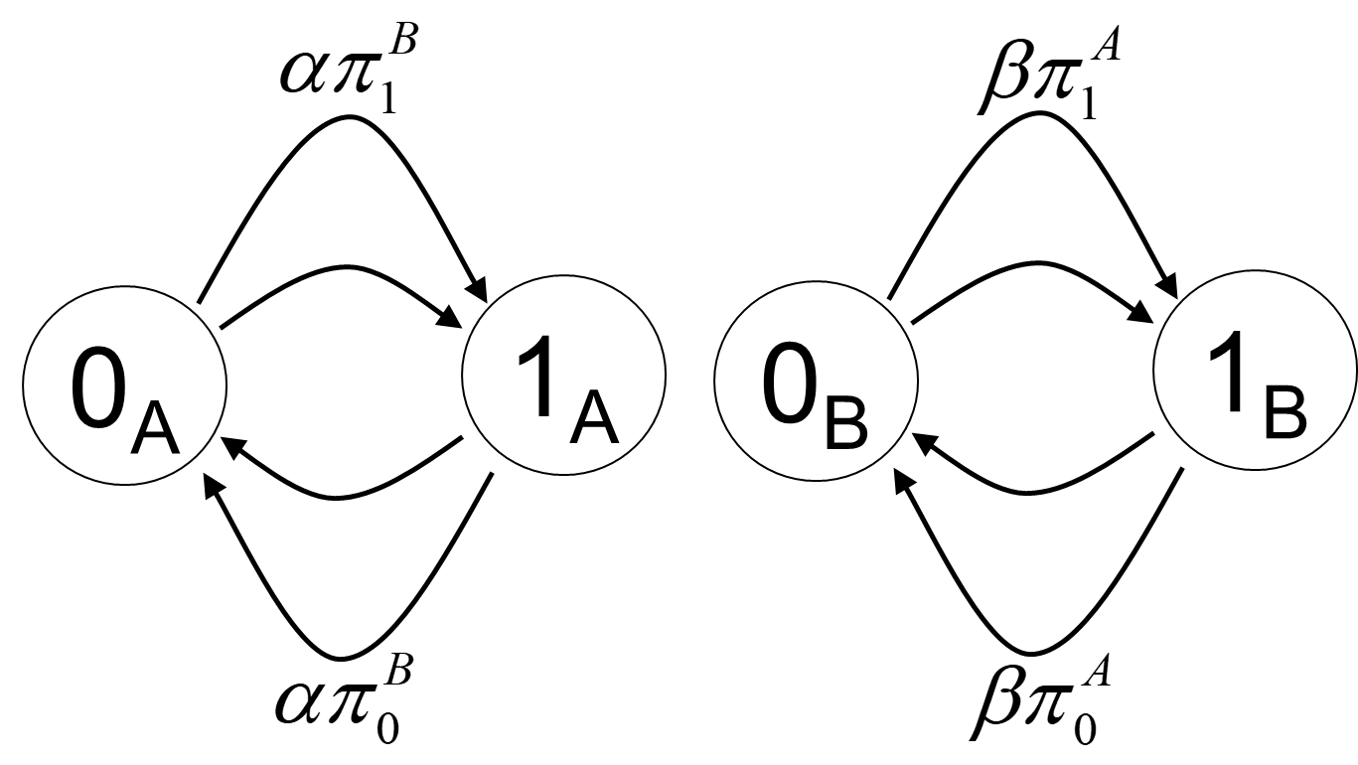}}
{Termination Scheme based on Partial Flow Conservation.\label{fig:fourstateMC3-termination}}
{}
\begin{flushleft}
\hspace{.5cm} \footnotesize \emph{Note}: Observe that $\pi^A$ and $\pi^B$ depend on each other.
\end{flushleft}
\end{figure}

\subsubsection{\textcolor{black}{Decomposition Procedure}}
\label{sec:decomposition}
\textcolor{black}{When decomposing the full MC, all decomposed subchains must maintain the steady-state distribution of the full MC; thus, we require $(\pi^A_0,\pi^A_1) \propto (\pi_{0A},\pi_{1A})$ and $(\pi^B_0,\pi^B_1) \propto (\pi_{0B},\pi_{1B})$, where ``$\propto$'' denotes that two vectors (distributions) are proportional. To satisfy this requirement, we \emph{terminate} (i.e., add extra transitions to) subchains after decomposition. Now, consider a termination scheme for $A$ (note: a symmetric argument holds for $B$) and let any added transitions be $\Delta^k_{i,j}$ from states $i_k$ to $j_k$, $k \in \{A,B\}$. A standard approach is to first check the return probability from each state in $A$ to each state in $A$ via states in $B$; this is accomplished by sample path analysis or an analysis of the transition matrix utilizing, for example, a censoring Markov chain technique. Once we find these return probabilities, we redirect transitions based on the return probabilities.}

\textcolor{black}{This scheme works fine, but the procedure tends to become complicated because a sample path depends on the structure of $B$. Our termination scheme does not require the use of return probabilities; rather, we impose a condition that each state in $A$ conserves \emph{partial} flow (i.e., net in/outflow), the condition necessary and sufficient to maintain the correct distribution of a subchain after decomposition. Note that the standard scheme using return probabilities satisfies this condition; however, there exist infinitely many alternative schemes that conserve partial flow and thus maintain the correct distribution. Our termination scheme is one of them; we redirect all outflows (from $A$) following the proportion of inflows (into $A$). In this simple model, the proportion of inflow into state $j_A$ is ${\pi}^B_j$ (extensions are straightforward for general cases). Thus, the added transitions in our termination scheme are $\Delta^A_{i,j}=\alpha {\pi}^B_j$ (see Figure \ref{fig:fourstateMC3-termination}). This termination scheme \emph{always} works regardless of the structure of $B$ because given we observe the distribution of $B$, our termination scheme is independent of transitions made within $B$. This is essentially the Markov property: the proportion of inflows to $A$ depends only on the \emph{present} distribution of $B$ and not on the \emph{past} transitions made within $B$. The use of the Markov property not only simplifies the termination scheme, but also reveal the impact of $B$ on $A$, as $\Delta^A_{i,j}$ is always represented as an explicit function of $\pi^B$ in our method.}

\subsubsection{\textcolor{black}{Aggregation Procedure}}
\label{sec:aggregation}
\textcolor{black}{Once we obtain $\pi^A$ and $\pi^B$, the next step is to derive the performance measure $E[f(X)]$. In a standard approach, we define ${p_A=\text{Pr}[X \in A]}$ and ${p_B=\text{Pr}[X \in B]}$ with a normalization condition $p_A+p_B=1$, find $p_A=\frac{\beta}{\alpha+\beta}$ and $p_B=\frac{\alpha}{\alpha+\beta}$ (for this model), retrieve the distribution of the full MC as $\pi=(p_A \pi^A_0,p_A \pi^A_1,p_B \pi^B_0,p_B \pi^B_1)$, and derive $E[f(X)]=\sum_{i \in \Omega} \pi_i f(i)$. This procedure is straightforward, but $\pi$, the normalized distribution of the full MC, tends to have a complicated expression that often makes the analytical representation of the performance measure of interest unwieldy. (Recall that our goal is to derive $E[f(X)]$, not $\pi$, unless required.)}

\textcolor{black}{In contrast, our method first finds the subchains' expectations $E_A[f(X)]$ and $E_B[f(X)]$, and then combines them using the total expectation theorem in MC settings (Theorem \ref{total expectation}), which in this example takes the form: $\frac{E[f(X)]}{p_A}=E_A[f(X)]+\frac{\alpha}{\beta}E_B[f(X)]$, where $p_A$ represents the reference state (set) probability, the probability of being in subchain $A$ in this example. This relationship between subchains and the original chain holds regardless of transition rates within each subchain; thus, we can focus on the analysis of each subchain to derive $E[f(X)]$. We do not need to impose a normalization condition separately, since it is part of this relationship (i.e., all expectations are 1 when ${f(X)=1}$), which is utilized to determine $p_A=(1+\frac{\alpha}{\beta})^{-1}$.}

\section{Total Expectation Theorem for an MC}
\label{sec:tet}
\subsection{Preliminaries}

We consider an ergodic (i.e., positive recurrent) continuous time MC that we decompose into multiple subchains, each of which is indexed with $j$, $j \in J^+$. We denote the set of states that compose a subchain $j$ as $A_j$ and the whole set of states in the full MC as $S$. These decomposed subchains $\{A_j: j \in J^+\}$ should be collectively exhaustive (i.e., a collection of decomposed subchains should form the full MC:  $S = \bigcup_{j \in J^+} {A_j}$), but not necessarily mutually exclusive (disjoint). In order to correct excess contribution from overlapping states in $\{A_j: j \in J^+\}$ when calculating a total expectation, we augment the set of subchains $\{A_j: j \in J^+\}$ using another set of subchains  $\{A_j: j \in J^-\}$ that satisfies the following condition, where $I(\cdot)$ represents the indicator function:
\begin{equation}
\label{Condition for a set J}
\sum_{j \in J^+} I(k\in A_j)-\sum_{j \in J^-} I(k\in A_j)=1, \forall k \in S.
\end{equation}
By considering all subchains in $J=J^+ \cup J^-$ (i.e., by adding up states in $\{A_j: j \in J^+\}$ and subtracting states in $\{A_j: j \in J^-\}$), we can make each state in $S$ contribute to the total expectation exactly once. \textcolor{black}{Under this framework, we allow some states to overlap more than once (e.g., in three or more subchains) as long as we subtract their contributions appropriately.} Note that if $\{A_j: j \in J^+\}$ forms a partition of the full set $S$, the index set $J^-$ becomes a null set; however, a set $J^+$ cannot be a null set. \textcolor{black}{As an example, if a level-dependent quasi-birth-death MC (LDQBD) changes birth and death rates at the threshold level $\ell$, then we may want to set $A_1$ and $A_2$ (both in $J^+$) as the sets of states up to the level $\ell$ and from the level $\ell$, respectively, and $A_3$ (in $J^-$) as a set of states with the level $\ell$. Alternatively, we can define $A_j$ to represent a subchain corresponding to the level (or phase) $j \in J^+=\{0,1,\cdots\}$ of the LDQBD queue, in which case $J^-=\emptyset$. In our method, we can define and use any set of subchains that satisfy Equation \eqref{Condition for a set J}.}

We denote stationary distributions of the full MC and a subchain $j$ as $\pi _k^{S}, \forall k \in S$, and $\pi _k^{A_j}, \forall k \in A_j$, respectively. We also denote expectations evaluated on the full MC and on a subchain $j$ as $E_S[\cdot]$ and $E_{A_j}[\cdot]$, respectively. In this paper we use the following simpler notations: $\pi _k \doteq \pi _k^S $, $\pi _k^j \doteq \pi _k^{A_j}$, $E[\cdot] \doteq E_S[\cdot] $, and $E_j[\cdot] \doteq E_{A_j}[\cdot]$.

Throughout the paper, we consider the case in which a decomposed subchain maintains the same stationary distribution as the full MC (up to a normalization constant): ${\pi _k} \propto \pi _k^j, \forall k \in {A_j}, \forall j \in J$. We call this condition \emph{conservation of distribution}. Obviously, a stationary distribution of a subchain is strongly affected by how we decompose the full MC and how we treat transition rates that are lost through the decomposition procedure. Our procedure that ensures conservation of the distribution of a decomposed subchain is discussed in \S\ref{sec:term}.

The conservation of distribution condition, or simply conservation of distribution, can be expressed in several different ways. To show certain conditions equivalent to conservation of distribution, we define the following quantities:

\begin{definition}
\label{ratio}
\it We define the following quantities that are functions of $\pi_k, \forall k \in S$:
\begin{enumerate}
  \item The ratio of steady-state probabilities: $ \beta_{kk'} \doteq \frac{\pi _k}{\pi _{k'}}.$
  \item The conditional steady-state probability (of the original MC) of being in a state $k$ given it is in a subchain $j$: $p_k^j \doteq \frac{\pi _k}{\sum\limits_{k' \in {A_j}} {\pi _{k'}}}.$
  \item The probability of being in a subchain $j$: $P_j \doteq \sum_{k \in {A_j}} {\pi _k}.$
  \item The conditional expectation of a function $f(\cdot)$ given that the MC is in a subchain $j$: $E[f(X)|{A_j}] \doteq \sum_{k \in {A_j}} {f(k)p_k^j}.$
\end{enumerate}
\end{definition}
Notice that these quantities are functions of ${\pi _k}$, not ${\pi _k^j}$. The conditions equivalent to conservation of distribution connect quantities defined by ${\pi _k}$ with functions of ${\pi _k^j}$, as we see in the following proposition:
\begin{proposition}
\label{conservation of distribution}
 Conservation of distribution: A stationary distribution is conserved for a subchain if and only if any of the following equivalent conditions are satisfied.
\begin{itemize}
  \item[] $(a) \qquad {\pi _k} \propto \pi _k^j,\forall k \in {A_j},\forall j \in J.$
  \item[] $(a') \qquad \beta_{kk'} = \frac{\pi _k^j}{\pi _{k'}^j} , \forall k \in A_j, \forall k' \in {A_j}, \forall j \in J.$
  \item[] $(b) \qquad p_k^j = \pi _k^j, \forall k \in {A_j}, \forall j \in J.$
  \item[] $(b') \qquad P_j = \frac{\pi _k}{\pi _k^j}, \forall k \in {A_j}, \forall j \in J.$
  \item[] $(c) \qquad E[f(X)|{A_j}] = E_j[f(X)], \forall j \in J \text{ and for any function } f(X).$
\end{itemize}
\end{proposition}

Condition $(c)$ in Proposition \ref{conservation of distribution} shows that under conservation of distribution, the change of probability measure does not affect the expected values, or in other words, a conditional expectation given the MC is in a specific subchain can be replaced by an expectation evaluated on a \emph{decomposed} subchain. This is an essential property when implementing the total expectation theorem in MC settings.

For notational convenience, under conservation of distribution, we define $\pi _{(k)}^j$, where a state $k$ does not necessarily belong to a subchain $j$:
$\pi _{(k)}^j \doteq \beta_{kk'} \cdot \pi_{k'}^j = \frac{\pi_k}{\pi_{k'}} \cdot \pi _{k'}^j, \forall k \in S, \forall k' \in A_j.$
Note that $\pi _{(k)}^j$ does not depend on $k'$ in subchain $j$. Note also that $\pi _{(k)}^j=\pi _k^j$ if $k \in A_j$. Using this notation, conditions $(a')$ and $(b')$ can be re-represented as follows:
\begin{itemize}
  \item[] $(a'') \qquad \beta_{kk'} = \frac{\pi _{(k)}^j}{\pi _{(k')}^j} , \forall k \in S, \forall k' \in S, \forall j \in J.$
  \item[] $(b'') \qquad P_j = \frac{\pi _k}{\pi _{(k)}^j}, \forall k \in S, \forall j \in J.$
\end{itemize}

\subsection{Total Expectation Theorem for an MC}
When we adopt a divide-and-conquer approach to an MC, we want a decomposed subchain to faithfully represent a part of the full MC; i.e., we require that a stationary distribution is conserved after decomposition: $\pi_k \propto \pi_k^j, \forall k \in A_j, \forall j \in J$. This conservation of distribution condition leads to the following theorem with help from Proposition \ref{conservation of distribution}. The theorem is essentially the total expectation theorem in MC settings.
\begin{theorem}
\label{total expectation}
Total expectation theorem for an MC: Under the conservation of distribution condition, an expectation of the full MC can be represented by expectations and steady-state probabilities of subchains, where a reference state $k$ can be any state in the full MC:
\begin{equation}
\label{tetequation}
\dfrac{E[ {f\left( X \right)}]}{\pi _k} =\sum\limits_{j \in J^+} {\dfrac{E_j [ {f(X)}]}{\pi _{(k)}^j}}  - \sum\limits_{j \in J^-} {\dfrac{E_j [ {f(X)}]}{\pi _{(k)}^j}}.
\end{equation}
\end{theorem}

\remark{\textcolor{black}{Equation \eqref{tetequation} holds for any $f(X)$. In particular, by setting $f(X)=1$, we obtain
\begin{equation*}
\label{identityequation}
\frac{1}{\pi_k} =\sum\limits_{j \in J^+} {\frac{1}{\pi _{(k)}^j}}  - \sum\limits_{j \in J^-} {\frac{1}{\pi _{(k)}^j}},
\end{equation*}
which determines $\pi_k$. This is an identity equation and serves as the normalization condition in the aggregation process of decomposed MCs.}}

\remark{Theorem \ref{total expectation} includes a reference state $k$, which can be arbitrarily chosen from the set $S$ but must be fixed throughout the calculation. This reference state can be replaced by a reference ``set'' $K$, which includes multiple states. The extension from a reference state $k$ to a reference set $K$ in Theorem \ref{total expectation} is straightforward by observing that $\pi_k$ and $\pi_k^j$ in Proposition \ref{conservation of distribution} can be replaced by $\pi_K(\doteq \sum_{k \in K} \pi_k)$ and $\pi_K^j(\doteq \sum_{k \in K \subseteq {A_j}} \pi_k^j)$, respectively.}

\textcolor{black}{Theorem \ref{total expectation} is the key to our decomposition method: Using Equation \eqref{tetequation}, we can combine all scaled expectations $E_j [{f(X)}]/\pi _{(k)}^j$ to derive the expectation of interest, and also obtain a normalization condition by setting $f(X)=1$. These desirable characteristics originate from the properties of the \emph{measure} (see Proposition \ref{additivity property}). In fact, by considering the $\sigma$-algebra $\mathcal{F}$ generated by a set $\{A_j, j \in J^+\}$, Proposition \ref{additivity property} immediately follows from Theorem \ref{total expectation}:}

\begin{proposition}
\label{additivity property}
Additivity property: Assume that conservation of distribution holds after a decomposition procedure. Let $\{A_j, j \geq 1\}$ be a collection of disjoint subchains (i.e., $J^+=\{j|j \geq 1\}$ and $J^-=\emptyset$) in the $\sigma$-algebra $\mathcal{F}$. Then a function $\mu_k(A_j)\doteq {E_j}[ {f(X)}]/{\pi _{(k)}^j}$ satisfies the countable additivity property:
$$\mu_k\left(\bigcup_{j=1}^{\infty} A_j \right) =\sum_{j=1}^{\infty} \mu_k(A_j).$$
\end{proposition}
This countable additivity property, together with assumptions, $\mu_k(\emptyset)=0$ and $\mu_k(A) \geq 0, \forall A \in \mathcal{F}$, indicates that the function $\mu_k(\cdot)$ satisfies the definition of a \emph{measure} (see, for example, \citealt*{halmos1950measure} and \citealt*{grimmett2001probability}). \textcolor{black}{Hence, $\mu_k(\cdot)$ inherits the following well-known properties that \textcolor{black}{all \emph{measure} possesses}.}
\begin{proposition}
\label{linear property}
\textcolor{black}{Let $A_j$ be a subchain of an MC. Then $\mu_k(A_j)\doteq {E_j}[ {f(X)}]/{\pi _{(k)}^j}$ satisfies the following properties:}

(1) For any subchain ${A_j}$  and its complement ${A_j}^c = S\backslash {A_j}$, $\mu_k(S)=\mu_k({A_j}) + \mu_k({A_j}^c).$

(2) If subchains ${A_2} \supseteq {A_1}$, then $\mu_k({A_2}) = \mu_k({A_1}) +\mu_k({A_2}\backslash {A_1}) \geq \mu_k(A_1).$

(3) For any two subchains ${A_1}$ and ${A_2}$, $\mu_k({A_1} \cup {A_2}) = \mu_k({A_1}) + \mu_k({A_2}) - \mu_k({A_1} \cap {A_2})$. In particular, if both subchains are disjoint, then  $\mu_k({A_1} \cup {A_2}) = \mu_k({A_1}) + \mu_k({A_2})$. 
\end{proposition}
\textcolor{black}{Propositions \ref{additivity property} and \ref{linear property} guarantee that any set of (e.g., disjoint or overlapping) subchains that is convenient for our analysis can be used, and the final result is obtained simply by summing up all parts and subtracting contributions from overlapping parts.}

\section{Termination}
\label{sec:term}
If we arbitrarily decompose an MC into subchains, each subchain independently analyzed is likely to have a different stationary distribution from its portion of the stationary distribution of the full MC: Proposition \ref{conservation of distribution} may not hold. This is because a subchain extracted from the full MC loses transitions between states in the subchain and states outside of the subchain. Such a loss can affect the distribution of a decomposed subchain; specifically, it could violate the conservation of distribution property. Therefore, in order to conserve the stationary distribution of a decomposed subchain, we must alter transitions at states that lose transitions. We call our alteration procedure \emph{termination}. \textcolor{black}{Countless such termination schemes exist, including the scheme obtained based on return probabilities by censoring (see the discussion in \S\ref{sec:decomposition}).} The simplest example of termination is \emph{truncation}, in which no alteration needs to be made for states that lose transitions. Truncation does not always conserve the stationary distribution of a decomposed subchain, unless the MC has a special structure (e.g., quasi-reversibility). In contrast, termination, if appropriately chosen, can \emph{always} conserve a stationary distribution. \textcolor{black}{In this section we demonstrate a procedure to obtain one of many possible termination schemes that satisfies both conservation of distribution and ergodicity.}

\subsection{Necessary and Sufficient Condition for a Decomposed Subchain to Conserve a Stationary Distribution}

We continue our analysis of the set $\{A_j$, $j \in J\}$, $J=J^+ \cup J^-$, which satisfies ${S}= \bigcup_{j \in J^+} {A_j}$. We denote a transition rate from state $k$ to state $k'$ in the full MC as ${q_{k,k'}}$. Let the set of states in $A_j$ that lose transitions due to decomposition be denoted as the boundary set: $\partial A_j= \left\{ {k \in {A_j}|{q_{k,k'}} + {q_{k',k}} > 0,\exists k' \in {A_j}^c} \right\}$. Let the set of states in $A_j$ that is not in $\partial A_j$ be denoted as the interior set: $\text{int}(A_j)(=A_j \backslash \partial A_j= \left\{ {k \in {A_j}|{q_{k,k'}} + {q_{k',k}} = 0,\forall k' \in {A_j}^c} \right\})$. Hence, $\partial {A_j} \cup {\mathop{\rm int}} ({A_j}) = {A_j}$ and $\partial {A_j} \cap {\mathop{\rm int}} ({A_j}) = \emptyset$. (Note that $\partial A_j \neq \emptyset$ since a subset $A_j (\subsetneq S)$ should communicate with its complement $A_j^c (=S\backslash A_j)$ under our ergodicity assumption for the full MC.)

We consider \emph{terminating} boundary states of a decomposed subchain: Adding new transitions $\Delta q_{k,k'}^j$ for some $k,k' \in \partial {A_j}$ in subchain $j$. (If $\Delta q_{k,k'}^j=0, \forall k \in \partial {A_j}, \forall k' \in \partial {A_j}$, then we call such a termination scheme truncation.) Termination of a decomposed subchain $j$ is added only to the states in $\partial {A_j}$ and not to any states in $\text{int}(A_j)$.

Under a termination scheme, we claim that the conservation of distribution condition shown in Proposition \ref{conservation of distribution} is equivalent to the \emph{boundary condition}, as described in the following proposition:

\begin{proposition} 
\label{boundary condition}
Boundary Condition: Under a termination scheme, conservation of distribution (Proposition \ref{conservation of distribution}) holds if and only if any of the following conditions are satisfied.

\begin{itemize}
\item[] $(d) \qquad {\pi _k} \propto \pi _k^j,\forall k \in {\partial A_j},\forall j \in J.$
\item[] $(d') \qquad  \beta_{kk'} = \frac{\pi _k^j}{\pi _{k'}^j} , \forall k \in \partial A_j, \forall k' \in {\partial A_j}, \forall j \in J.$
\item[] $(e) \qquad p_k^j = \pi _k^j, \forall k \in {\partial A_j}, \forall j \in J.$
\item[] $(e') \qquad P_j = \frac{\pi _k}{\pi _k^j}, \forall k \in {\partial A_j}, \forall j \in J.$
\end{itemize}
\end{proposition}

According to Proposition \ref{boundary condition}, the condition for conservation of distribution reduces to the boundary condition:
${\pi _k} \propto \pi _k^j,\forall k \in {A_j} \Leftrightarrow {\pi _k} \propto \pi _k^j,\forall k \in \partial {A_j}.$
This is intuitively obvious because termination 
only changes flow balance equations at boundary states. As long as the boundary distribution is maintained, the entire distribution of a decomposed subchain is conserved. In this paper, we use the terms ``the conservation of distribution condition'' and ``the boundary condition'' interchangeably as they are equivalent.

How can we satisfy the boundary condition? Needless to say, the simplest way to satisfy this condition is to impose $\pi _k^j$ follow ${\pi _k} \propto \pi _k^j,\forall k \in \partial {A_j}$, if we know how ${\pi _k}$ is distributed at the boundary states. This is often possible when these boundary states also belong to another subchain, for which we know the distribution \textcolor{black}{\citep[as an example, see][]{doroudi2017clearing}}. However, in general, we need to find $\pi _k^j$ at boundary states by solving an independent subchain $j$ with termination. We next discuss the condition for termination that ensures the boundary condition.

\subsection{Partial Flow Conservation: An Equivalent Condition to the Boundary Condition}
One obvious conclusion we can draw from Proposition \ref{boundary condition} is that if the number of states in a boundary set is one for a subchain (i.e., $|\partial {A_j}|=1, j \in J$), then the boundary condition for the subchain is automatically satisfied, and therefore, the stationary distribution of this subchain is conserved regardless of any termination schemes we apply, including the simplest termination scheme: Truncation. Intuitively, if a subchain connects to the rest of the chain through a single state, the flow out from the subchain through the single state must be balanced by the flow into the subchain through the same single state, in which case a loss of flows due to decomposition does not alter the distribution of a decomposed subchain.

\textcolor{black}{However, if there is more than one boundary state, we need to control flows at each boundary state to potentially compensate for a loss of transitions due to decomposition. Specifically, in order to conserve the boundary distribution (and thus the original stationary distribution up to a normalization constant) after decomposition, we may need to add extra transitions $\Delta q^j_{k,k'}$ (\emph{termination}) to conserve the net in/outflow (partial flow) at each boundary state. The following lemma provides the necessary and sufficient condition for termination to satisfy:}
\begin{lemma}
\label{partial flow}
Partial flow conservation: The stationary distribution is conserved after decomposition if and only if termination $\Delta q^j_{k,k'}$ conserves partial flow at all boundary states:

$(f) \qquad {\pi _k}\sum\limits_{k' \in A_j^c} {{q_{k,k'}}}  - \sum\limits_{k' \in A_j^c} {{\pi _{k'}}{q_{k',k}}}  = {\pi _k}\sum\limits_{k' \in \partial {A_j}} {\Delta {q_{k,k'}^j}}  - \sum\limits_{k' \in \partial {A_j}} {{\pi _{k'}}\Delta {q_{k',k}^j}}, \forall k \in \partial {A_j}, \forall j \in J.$
\end{lemma}\smallskip
\textcolor{black}{Lemma \ref{partial flow} shows that if termination $\Delta q_{k,k'}^j$ for a subchain $j$ can replicate the net outflow to (or net inflow into) the complementary set at all boundary states of the subchain $j$, then the decomposed subchain $j$ with added transitions $\Delta q_{k,k'}^j$ conserves the stationary distribution of the full MC. For example, termination $\Delta^A_{k,k'}=\alpha\pi^B_{k'}$ and $\Delta^B_{k,k'}=\beta\pi^A_{k'}$ for the four-state CTMC in \S\ref{sec:four-state} (and also $\Delta^A_{k,k'}=\alpha_k \hat{\pi}^B_{k'}$ and $\Delta^B_{k,k'}=\beta_k \hat{\pi}^A_{k'}$ for the four-state CTMC in Appendix \ref{sec: extension}) are selected so that they satisfy the partial flow conservation condition (Lemma \ref{partial flow}) and thus maintain the stationary distribution of the full MC after decomposition.}

\textcolor{black}{The remaining task is to develop a scheme to find a termination that satisfies the condition $(f)$. Before discussing a general case, consider a special case, partial flow balance, which was first introduced by  \cite{whittle1968equilibrium} and was discussed extensively in the context of queueing networks by \cite{kellyreversibility}. Assuming all states satisfy the partial flow balance condition, ${\pi _k}\sum\limits_{k' \in A_j^c} {{q_{k,k'}}}  - \sum\limits_{k' \in A_j^c} {{\pi _{k'}}{q_{k',k}}}  = 0,\forall k \in \partial {A_j},$ it is obvious that a simple truncation (i.e., $\Delta {q_{k',k}^j}=0, \forall k, k' \in \partial A_j$) is sufficient to satisfy Lemma \ref{partial flow} (and thus conserve the distribution). This property is known as the \emph{state truncation property} (see \citealt{nelson1995probability}). (Note: If ergodicity is not maintained by truncation, we still need termination.) However, in general, partial flow balance does not always hold and/or ergodicity may not be maintained by truncation. In such cases, we need to find a non-trivial termination that satisfies the partial flow conservation condition $(f)$.}

\subsection{Termination That Satisfies Partial Flow Conservation}
\label{sec:termination}
\textcolor{black}{In this subsection, we discuss our scheme to obtain a termination that satisfies the partial flow conservation condition $(f)$ in Lemma \ref{partial flow}, the necessary and sufficient condition to conserve the distribution after decomposition. Recall that our goal is to find one of many possible termination schemes that satisfy $(f)$. For this purpose, we first show the sufficient condition to satisfy $(f)$ (Proposition \ref{termination}) and then show one of the termination schemes that satisfy this sufficient condition (Corollary \ref{special termination}). We continue to require termination to maintain ergodicity of a decomposed subchain.}

\begin{proposition}
\label{termination}
Termination condition sufficient to satisfy Lemma \ref{partial flow}: The boundary distribution is conserved if termination $\Delta {q_{k,k'}^j}$ satisfies the following flow conservation conditions at all boundary states:

(1) Outflow (from a subchain) condition: $\sum\limits_{k' \in A_j^c} q_{k,k'} = \sum\limits_{k' \in \partial {A_j}} {\Delta {q_{k,k'}^j}}, \forall k \in \partial {A_j}.$

(2) Inflow (to a subchain) condition: $\sum\limits_{k' \in A_j^c} {\pi _{k'}}{q_{k',k}}  = \sum\limits_{k' \in \partial {A_j}} {\pi _{k'}}\Delta {q_{k',k}^j}, \forall k \in \partial {A_j}.$
\end{proposition}

\remark{\textcolor{black}{Since Proposition \ref{termination} is a sufficient condition (for Lemma \ref{partial flow}), termination satisfying Lemma \ref{partial flow} may not satisfy Proposition \ref{termination}.} For example, we know we can arbitrarily add or drop self-transitions from termination without affecting Lemma \ref{partial flow}. However, such a change in termination is not allowed by Proposition \ref{termination}. Another example is that if we could find another termination $\Delta {q'}_{k,k'}^j$ that satisfies
$$\pi_k \sum\limits_{k' \in \partial {A_j}} {\Delta {{q'}_{k,k'}^j}}= \sum\limits_{k' \in \partial {A_j}} {{\pi _{k'}}\Delta {{q'}_{k',k}^j}}, \forall k \in \partial A_j,$$ 
then the new termination $\Delta {q_{k,k'}^j}+\Delta {{q'}_{k,k'}^j}$ does not satisfy Proposition \ref{termination} but satisfies Lemma \ref{partial flow}.}

\textcolor{black}{To describe our termination scheme, we define two key parameters: (i) the total outflow rate from state $k \in \partial A_j$ to $A_j^c$: $q^j_k \doteq  \sum\limits_{m \in A_j^c} q_{k,m}$, and (ii) the proportion of inflow (the \emph{normalized} relative frequency of visits) to state $k' \in \partial A_j$ from $A_j^c$: $\hat{\pi}^{-j}_{k'} \doteq \dfrac{\sum_{m \in A_j^c} {\pi_m^{-j}}{q_{m,k'}}}{\sum_{k'' \in \partial A_j}\sum_{m \in A_j^c} {\pi_m^{-j}}{q_{m,k''}}}$. Here, we assume that we know the correct distribution $\pi _{m}^{-j}$ of properly terminated $A_j^c$, which maintains ${\pi _m^{-j} \propto \pi_m, \forall m \in A_j^c}$. Since only the boundary states of $A_j^c$ (that have non-zero transition rates to subchain $j$) are considered in the summation and $\pi_m^{-j}$ appears both in the numerator and denominator in the expression, we could alternatively replace $\pi_m^{-j}$ with the boundary-state distribution of $A_j^c$. In this paper we use the same expression $\pi _{m}^{-j}$ for both the full and boundary-state distributions of properly terminated $A_j^c$ unless it creates confusion; we use the convenient alternative depending on a model. The following corollary shows one possible termination scheme:
\begin{corollary}
\label{special termination}
Termination scheme: The following termination $\Delta {q_{k,k'}^j}$ satisfies Proposition \ref{termination} (and thus Lemma \ref{partial flow}):
\begin{equation}
\label{general termination}
\Delta {q_{k,k'}^j}=q^j_k \hat{\pi}^{-j}_{k'}, \forall k, k' \in \partial A_j.
\end{equation}
\end{corollary}}
\remark{Note again that the termination scheme in Corollary \ref{special termination} can be modified to include or exclude self-transitions because such a modification does not affect Lemma \ref{partial flow} \textcolor{black}{(although it may not satisfy Proposition \ref{termination})}.}
\remark{If we can partition the external states (states in $A_j^c$) into multiple non-communicating classes, Corollary \ref{special termination} can be applied to each class, and we can obtain termination schemes for multiple classes. Aggregation of these terminations represents the correct termination that satisfies Lemma~\ref{partial flow}.}

\textcolor{black}{Corollary \ref{special termination} simply states that one possible scheme to satisfy Proposition \ref{termination} and thus Lemma \ref{partial flow} is to redirect the total outflow from state $k \in \partial {A_j}$ to state $k' \in \partial {A_j}$ proportional to the inflow to state $k' \in \partial {A_j}$. It is intuitively obvious that the outflow condition (1) of Proposition \ref{termination} is satisfied because we redirect all outflows from $k \in \partial {A_j}$, and that the inflow condition (2) is satisfied because we redirect all outflows following the inflow proportion at state $k' \in \partial {A_j}$. For example, the four-state CTMC in \S\ref{sec:four-state} shows $q_k^A=\alpha$ and $\hat{\pi}^{-A}_{k'} =\pi_{k'B}/(\pi_{0B}+\pi_{1B}) = \pi^B_{k'}$ for $k,k' \in \{0,1\}$, given ${\pi^B_{k'} \propto \pi_{k'B}}, k' \in \{0,1\}$. Thus, to maintain the distribution of $A$ after decomposition, we can set termination ${\Delta^A_{k,k'}=q_k^A \hat{\pi}^{B}_{k'} =\alpha \pi^B_{k'}}$ for $k,k' \in \{0,1\}$, assuming $\pi^B_{k'}$ is the correct distribution of $B$. (Another four-state CTMC example is presented in Appendix \ref{sec: extension}.)}

\textcolor{black}{The termination scheme following Corollary \ref{special termination} has two advantages. First, finding the termination $\Delta {q_{k,k'}^j}$ is a simple task because $\Delta {q_{k,k'}^j}$ is the product of two functions, each of which depends only on a single parameter (either $k$ or $k'$). Second, the termination $\Delta {q_{k,k'}^j}$ explicitly depends only on $\pi_m^{-j}$ and not on the structure of $A_j^c$. Thanks to the second property and the Markov property, we can always construct $\Delta {q_{k,k'}^j}$ properly given $\pi_m^{-j}$.}

\textcolor{black}{For problems that are analytically solvable, the proportion $\hat{\pi}^{-j}_{k'}$ takes a simpler form. For example, if a single state $s$ in subchain $j$ has the only inflow from $A_j^c$, then $\hat{\pi}^{-j}_{k'}=\delta_{k' s}$ (Kronecker delta), and thus termination is reduced to a simple redirection of all outflows into state $s$: ${\Delta {q_{k,k'}^j} = q^j_k \delta_{k' s},} {\forall k,k' \in \partial{A_j}}$, which does not require any knowledge of $\pi_m^{-j}$; we see this example in \S\ref{sec: CBS}. Another example is when the proportion $\hat{\pi}^{-j}_{k'}$ is a simple function of $\pi_m^{-j}$, which we observe in the example in \S\ref{sec:mtmt1queue}. However, in general, an appropriate termination is a complex function of $\pi_m^{-j}$, and thus we may need to resort to a recursive method; for this purpose, we provide a numerical algorithm to computationally implement Corollary \ref{special termination} in \S\ref{sec:algorithm}.}

\section{Applications of Our MC Decomposition Method}
\label{sec:appl}
\textcolor{black}{In this section we show several applications of our method. First, we present our decomposition algorithm based on Corollary \ref{special termination} and apply the algorithm in numerical experiments. Next we analyze a queueing system with adjustable staffing levels, which has been used to model the border gates/toll booths system without customer abandonment (\citealt{zhang2009performance} and \citealt{bhandari2008exact}). Lastly, we analyze a single-server Markov-modulated MC, an M\textsubscript{t}/M\textsubscript{t}/1 queue. This model belongs to a larger class (\emph{class} $\mathbb{M}$) of skip-free (in level), unidirectional (in phase) quasi-birth-death (QBD) MCs (\citealt{doroudi2017clearing}). There are additional applications; interested readers can refer to, e.g., Appendix \ref{sec: relationship} and \cite{sasanuma2021approximate} which describes the analysis of a single station two-stage reneging queue.}
\textcolor{black}{\subsection{Decomposition Algorithm}
\label{sec:algorithm}}
\textcolor{black}{We formulate our termination scheme using matrix notation and present our decomposition algorithm. We consider decomposing an MC into disjoint sets, noting the following two cases.}
\textcolor{black}{\subsubsection{Case 1 (Flow balance holds at each cut)}
When the flow balance condition holds at each cut of the two subchains, a termination scheme is simplified. Consider decomposing a full MC into two subchains, in which case the condition must hold because the flow balance condition at the cut of the two subchains is reduced to the global flow balance condition. Let sets $A_{i} ,i\in \{ 1,2\} $ be a partition of a full MC, and $w^i$ be the probability (weight) that the current state is in each corresponding set. Denote the infinitesimal generator matrix (transition rate matrix) $Q$ of the original MC as $$Q=\left(\begin{array}{cc} {Q_{11}} & {Q_{12}} \\ {Q_{21}} & {Q_{22}} \end{array}\right),$$ where a subscript $i$ represents subchain $i$. We want to find a termination $\Delta Q_{ii} $ for $i=1, 2$. Define $E_{i} $ as a matrix of ones that has the same size as $Q_{ii}$. Denote the steady-state probability of decomposed subchain $i$ as a row vector $\pi ^{i} $. We introduce row vectors  $o^{i} $, $\xi^{i} $ and column vectors $e^{i} $, $q^{i} $, all of which have the same dimension as $\pi ^{i}$. Each element of $o$, $e$, $\xi$, and $q$ are defined as follows: $o_{k}^{i} =0$ and $e_{k}^{i} =1$ for all states $k$ in chain $i$. $\xi_{k}^{i} $ is an inflow proportion, a relative frequency of visits to state $k$ in chain $i$ from outside (of chain $i$), satisfying $\xi^{1} =\dfrac{\pi ^{2} Q_{21}}{\pi ^{2} Q_{21} e^{1}} , \; \xi^{2} =\dfrac{\pi ^{1} Q_{12} }{\pi ^{1} Q_{12} e^{2} }.$
$q_{k}^{i} $ is a total transition rate out from state $k$ in chain $i$ to outside (of chain $i$), satisfying $q^{1} =Q_{12} e^{2} , \; q^{2}=Q_{21} e^{1}.$}

\textcolor{black}{Following Corollary \ref{special termination}, with these notations an appropriate termination can be represented as
\begin{equation}
\label{deltaQ}
\Delta Q_{11} =q^{1}\xi^{1} =Q_{12} e^{2}\frac{\pi ^{2} Q_{21} }{\pi ^{2} Q_{21}e^{1} }, \quad \Delta Q_{22} =q^{2}\xi^{2} =Q_{21} e^{1}\frac{\pi ^{1} Q_{12}}{\pi ^{1} Q_{12}e^{2}}.
\end{equation}}\textcolor{black}{Note that the terminated $Q$-matrix, $Q_{ii} + \Delta Q_{ii}$, remains a generator matrix since the sum of all elements of any row of $Q_{ii} + \Delta Q_{ii}$ is conserved at 0 (because, following our termination scheme, transitions to all complementary sets are redirected to the states in chain $i$). From the balance equation ${\pi ^{i} \left(Q_{ii} +\Delta Q_{ii} \right)=o^{i}}$ and the normalization condition $\pi ^{i} E_{i} =e^{i} $, we know that $\pi ^{i} $ should satisfy ${\pi ^{i} \left(Q_{ii} +\Delta Q_{ii} +E_{i} \right)=e^{i}}$, or ${\pi ^{i} =e^{i} \left(Q_{ii} +\Delta Q_{ii} +E_{i} \right)^{-1}}$. Therefore, we obtain two equations that should hold simultaneously:
\begin{equation}
\label{findpi}
\pi ^{1} =e^{1} \left(Q_{11} +\Delta Q_{11}+E_{1} \right)^{-1}, \quad \pi ^{2} =e^{2} \left(Q_{22} +\Delta Q_{22}+E_{2} \right)^{-1}.
\end{equation}}\textcolor{black}{To solve this set of equations, we start from an arbitrary distribution of $\pi ^{i} $ and find the correct $\pi ^{i}$ recursively. \textcolor{black}{Once all $\pi ^{i} $ have converged, we obtain their weights $w^i$ using the flow balance condition $w^1 \pi^1 Q_{12} e^2 = w^2 \pi^2 Q_{21} e^1$ with the normalization condition $w^1+w^2=1$, or using the total expectation theorem (Theorem \ref{total expectation}) $\frac{E[f(X)]}{w^1}=E_1[f(X)]+\frac{\pi^1 Q_{12} e^2}{\pi^2 Q_{21} e^1} E_2[f(X)]$ with $f(X)=1$:
\begin{equation}
\label{findw}
w^1=\frac{\pi^2 Q_{21} e^1}{\pi^1 Q_{12} e^2+\pi^2 Q_{21} e^1},  \quad w^2=\frac{\pi^1 Q_{12} e^2}{\pi^1 Q_{12} e^2+\pi^2 Q_{21} e^1}.
\end{equation}}
We summarize our algorithm to find the steady-state probabilities of the subchains as follows:\\
Step 1. Set the initial distributions for all subchains (e.g., a uniform distribution).\\
Step 2: Compute $\pi ^{i}$ using Equations \eqref{deltaQ} and \eqref{findpi} based on $\pi ^{j}, j\ne i$.\\
Step 3: Repeat Step 2 until all $\pi ^{i}$ converge (i.e., a stopping criterion is met).\\
\textcolor{black}{Step 4: Compute $w^i$ using Equation \eqref{findw} and obtain the solution $\pi=(w^1 \pi^1, w^2 \pi^2)$.}}
\textcolor{black}{\subsubsection{Case 2 (Flow balance does not hold at each cut)}
When the flow balance condition does not hold at a cut between two subchains, \textcolor{black}{which is often the case when we decompose a full MC into three or more subchains, a termination $\Delta Q_{ii}$ requires the knowledge of both the distribution $\pi^j$ and its probability weight $w^j$ of subchain $j$ for all $j \neq i$. To formulate this case, we consider an MC composed of three disjoint subchains 1, 2, and 3. (It is straightforward to extend this argument to a case with more than three subchains.) We use the same notations as the case (1). The relative frequency to visit the states in subchain 1 from subchains 2 and 3 (i.e., the proportion of inflow into subchain 1) is now}
\[\xi^{1} =\frac{w^{2}}{w^{2} +w^{3}} \cdot \frac{\pi ^{2} Q_{21}}{\pi ^{2} Q_{21}e^{1}} +\frac{w^{3}}{w^{2} +w^{3}} \cdot \frac{\pi ^{3} Q_{31}}{\pi ^{3} Q_{31}e^{1} }.\]
Or by defining $w_{-1}^{2} =\dfrac{w^{2} }{w^{2} +w^{3}}, w_{-1}^{3} =\dfrac{w^{3} }{w^{2} +w^{3}}$, we have a simpler form:
\[\xi^{1} =\frac{w_{-1}^{2} \pi ^{2} Q_{21}}{\pi ^{2} Q_{21}e^{1}} +\frac{w_{-1}^{3} \pi ^{3} Q_{31}}{\pi ^{3} Q_{31}e^{1}}.\]
Correspondingly, the total transition rate (outflow) to subchains 2 and 3 from subchain 1 becomes
\[\left(q^{1} \right)'=Q_{12} e^{2} +Q_{13} e^{3}.\]
Note that this total transition rate can be calculated just once at the beginning of each iteration.
Following Corollary \ref{special termination}, an appropriate termination is
$$\Delta Q_{11} =\left(q^{1} \right)'\xi^{1} =\left(Q_{12} e^{2} +Q_{13} e^{3}\right)\left(\frac{w_{-1}^{2} \pi ^{2} Q_{21} }{\pi ^{2} Q_{21}e^{1}} +\frac{w_{-1}^{3} \pi ^{3} Q_{31}}{\pi ^{3} Q_{31}e^{1}} \right).$$
Thus, the equation for $\pi ^{i} $ becomes
\begin{equation}
\label{deltaQ2}
\pi ^{i} =e^{i} \left(Q_{ii} + \Delta Q_{ii}+E_{i} \right)^{-1}, \text{where } \;
\Delta Q_{ii}=\left(\sum _{j\ne i}Q_{ij} e^{j} \right)\left(\sum _{k\ne i}\frac{w_{-i}^{k} \pi ^{k} Q_{ki} }{\pi ^{k} Q_{ki}e^{i}}\right).
\end{equation}}

\textcolor{black}{After we update $\pi ^{1} $, we need to update $w^{1}$. This $w^{1}$ is obtained using the global balance condition (i.e., the total flow out from subchain 1 is equal to the total flow into subchain 1), which is represented as
$$w^{1} \left(\pi ^{1} Q_{12}e^{2} +\pi ^{1} Q_{13}e^{3} \right)=w^{2}\pi ^{2} Q_{21}e^{1} +w^{3}\pi ^{3} Q_{31}e^{1},$$ or
\begin{equation}
\label{weight}
w^{i} =\frac{\left(\sum _{k\ne i}w^{k} \pi ^{k} Q_{ki} \right)e^{i} }{\pi ^{i} \left(\sum _{j\ne i}Q_{ij} e^{j}  \right)}.
\end{equation}}

\textcolor{black}{We obtain a set of equations, Equations \eqref{deltaQ2} and \eqref{weight}, with unknown $\pi^i$ and $w^i$ that should hold simultaneously. To solve this set of equations, we start from an arbitrary distribution of $\pi ^{i}$ and weight $w^i$, and find correct $\pi ^{i}$ and $w^i$ recursively. Once all $\pi ^{i}$ have converged, which will imply that all $w^i$ have also converged, we normalize $w^i$ (since Equation \eqref{weight} determines only the proportion of $w^i$) and obtain the solution $\pi$. We summarize our algorithm as follows.\\
Step 1. Set the initial distributions and weights for all subchains.\\
Step 2: Compute $\pi ^{i}$ and $w^i$ using Equations \eqref{deltaQ2} and \eqref{weight} based on $\pi ^{j}, w^j, j\ne i$.\\
Step 3: Repeat Step 2 until all $\pi ^{i}$ converge (i.e., a stopping criterion is met).\\
Step 4: Normalize $w^i$ and obtain the solution using $\pi=(w^1 \pi^1, w^2 \pi^2, \cdots)$.}
\textcolor{black}{\subsubsection{Computational Experiments}
The two key steps in our algorithm, a computation of $w^i$ using Equation \eqref{weight} and that of $\pi^i$ using Equation \eqref{deltaQ2}, correspond to aggregation and disaggregation procedures, respectively, in an iterative aggregation/disaggregation (IAD) method. Thus, our algorithm belongs to the family of IAD algorithms, whose convergence to the stationary solution is known to be geometric (exponential) (\S6 in \citealt{stewart1994introduction}). We confirmed our algorithm's exponential convergence to a stationary solution for all $Q$-matrices we tested: e.g., small/large, sparse/dense, nearly completely decomposable (NCD), and others.}
\begin{figure}[h]
\FIGURE
{\includegraphics*[scale=.9]{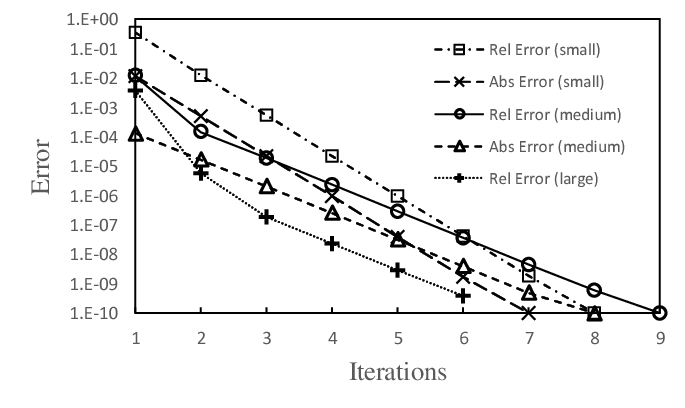}}
{\textcolor{black}{Relative error and absolute error of our algorithm depending on the number of iterations.}\label{fig:iterations}}
{}
\begin{flushleft}
\hspace{-.3cm} \textcolor{black}{\footnotesize \emph{Notes}: Sizes of $Q$-matrix are small ($n=8$), medium ($3,000$), and large ($30,000$). The precision of our numerical results is 1.E-10. We did not evaluate $E_{abs}$ for $n=30,000$ because $\pi$ was unavailable following a direct approach.}
\end{flushleft}
\end{figure}

\textcolor{black}{Here we show selected results of our experiments. We coded the algorithm in Python and tested it for the three subchain case with $Q$-matrices ($n\times n$ matrix) whose size $n \in \{8, 3000, 30000\}$. The smallest is the Courtois model that is frequently used to test operations of algorithms that solve MCs. The $Q$-matrix of the Courtois model is $8 \times 8$ and shows an NCD structure, which we split into three decomposed subchains with \{3,2,3\} states, respectively. For other $Q$-matrices, each non-diagonal element of the matrices is randomly chosen between 0 and 1 (thus, the matrices are very dense) and each diagonal term is set to make the sum of row elements equal to 0; we split the corresponding MC into three subchains with an equal number of states.}

\textcolor{black}{We evaluate our algorithm's performance by measuring errors (in $L_1$ norm) of the full MC's steady-state probability distribution. Specifically, at the end of each iteration, we obtain the full distribution $\pi_{(k)}$ (at the $k\textsuperscript{th}$ iteration) and compare it with the true full distribution $\pi$, the correct solution obtained using the direct approach within the computer's floating point precision. (Note that, in practice, it is not necessary to obtain the approximation $\pi_{(k)}$ at each iteration or the true distribution $\pi$ since all we need is the convergence of all subchains' distributions $\pi^i$. In this experiment, however, we obtained $\pi_{(k)}$ and $\pi$ because we need to evaluate the performance of our algorithm.) We use the following two indices:}

\textcolor{black}{$\bullet$ Relative error $E_{rel} \doteq ||\pi_{(k-1)}-\pi_{(k)}||_1$}

\textcolor{black}{$\bullet$  Absolute error $E_{abs} \doteq ||\pi-\pi_{(k)}||_1$}\\
\textcolor{black}{Numerical results are shown in Figure \ref{fig:iterations}, which confirms that both relative and absolute errors converge to zero with an (almost) constant exponential rate as we increase the number of iterations.}

\textcolor{black}{There are many possible extensions of our algorithm: For example, Step 2 can be computed in parallel to speed up the execution of our algorithm; each $\pi^i$ may be computed iteratively using, for example, the Gauss-Seidel method, or evaluated by further decomposing a subchain. In fact, many implementation details explained in \cite{stewart1994introduction} are also applicable to our algorithm. Such considerations will be an interesting topic for future research.}

\subsection{Queueing System with Adjustable Staffing Levels: CBS Model}
\label{sec: CBS}
We analyze an MC representing a queueing system with adjustable staffing levels. This queueing system is used to analyze the performance of border gates/toll booths/server farms. The staffing level is controllable, but operators may not want to adjust staffing levels frequently because changing staffing levels is costly. The simplest control policy for adjusting staffing levels is called the Congestion-Based Staffing (CBS) policy \citep*{zhang2009performance}. Under this CBS policy, if the queue length is decreased, some toll gates are closed, and if increased, some gates are opened. (We assume that each gate is operated by one staff member.) The most basic CBS policy uses two thresholds: The lower threshold $n$, which switches the system to operate with a lower number of staff, and the upper threshold $N$, which switches the system to operate with a larger number of staff. Operators of toll gates seek the optimal $(n,N)$ combination to minimize the sum of three costs: The cost of staffing, customers' waiting time, and switching staffing levels. The CBS policy is commonly used in practice, but its exact solution presented in \cite{zhang2009performance} is very complicated and hence an approximation is often sought. As an example to illustrate the simplicity of our MC decomposition method and its capability to cope with more extended models, we apply our method to derive simple, exact, analytical expressions for operational costs in the CBS model.

We assume there is no customer abandonment (neither reneging nor balking). The total number of toll gates are $c$, where $e$ out of $c$ ($c>e>0$) are extra gates that can be opened or closed (all at the same time). Arrivals of vehicles to the toll gate occur at constant rate $\lambda$, according to a Poisson process. Service times at each gate are distributed exponentially with parameter $\mu$. Under the CBS policy, $e$ gates are opened when the total number of vehicles in the system reaches the upper threshold $N$ and are closed when the total number of vehicles in the system reaches the lower threshold $n$, where $N>n$ must hold. We limit our analysis to the case where $n \geq c$ and arrival rate is constant. Extensions of this model are discussed in \cite{zhang2009performance} (allowing the lower threshold $n$ to be less than $c$) and \cite{bhandari2008exact} (time-varying arrival rate). The $n<c$ case in the CBS model can be handled using the same approach we present here. Analysis of the model with time-varying arrival rate requires a numerical decomposition approach, where we recursively identify a combination of $n$ and $N$ that optimizes the total operational cost. This analysis will be a potential research topic in the future.

\begin{figure}[h]
\FIGURE
{\includegraphics*[scale=0.4]{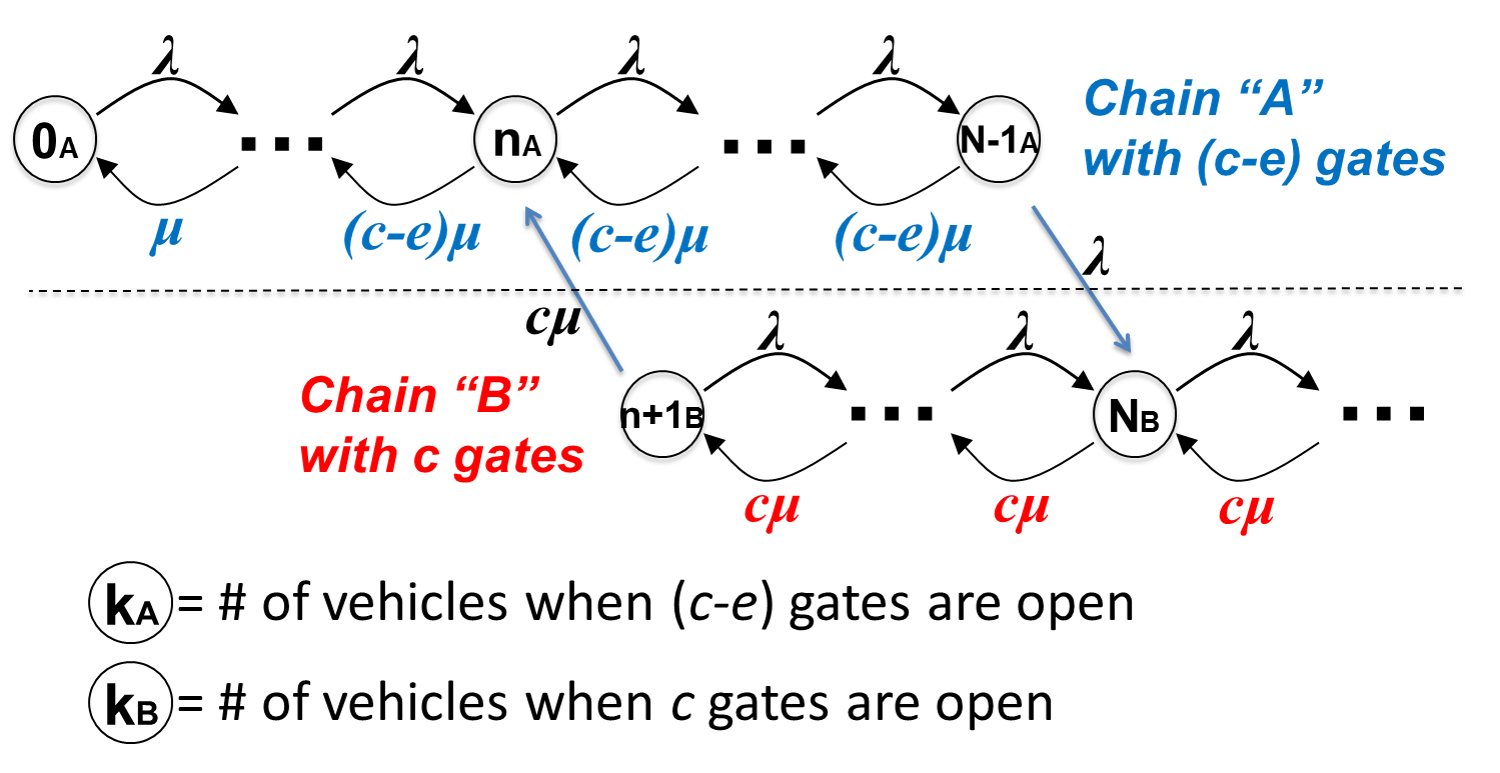}}
{The CBS model.\label{fig:MarkovChain-CBS}}
{}
\end{figure}

Figure \ref{fig:MarkovChain-CBS} shows the MC model corresponding to a toll gate system with a CBS policy. We call the upper part of the MC in Figure \ref{fig:MarkovChain-CBS} chain A and the lower part chain B. Notice that chains A and B represent M/M/c-e/N-1 and M/M/1 queues, and correspond to low capacity and high capacity modes, respectively. As Figure \ref{fig:MarkovChain-CBS} shows, a system switches to high capacity mode (chain B) when an arrival occurs at state $N-1$ in chain A (which we denote $N-1_A$), and switches to low capacity mode when a departure of a vehicle occurs at state $n+1$ in chain B ($n+1_B$). A \emph{changeover cycle time} corresponds to the time between two consecutive departure times from Chain A to Chain B. Define $c_g$, $c_w$, and $c_s$ as the cost of operation per gate per unit time, the waiting cost per vehicle per unit time, and a switching cost per cycle, respectively. Define $P_A$, $P_B(=1-P_A)$, $L_Q$, $f$, and $T_{cycle}$ as the probability to be in chain A, the probability to be in chain B, the expected number of vehicles in a queue, the frequency of changeover, and one changeover cycle time. The total cost is the sum of the following costs:

1. Average staffing cost: $C_{staff} = c_g \cdot ( c-e +  e \cdot P_B),$

2. Average waiting cost: $C_{wait} =c_w \cdot L_Q,$

3. Average switching cost: $C_{switch} =c_s \cdot E[f]=c_s \cdot \frac{1}{E[T_{cycle} ]} =c_s \cdot \lambda \cdot \pi_{N-1_{A}}.$\\
Note that $1/{E[T_{cycle}]} =\lambda \pi_{N-1_{A}}$ is used above. This is derived by applying Little's law to the ergodic closed system with 1 job. 

\subsubsection{Analysis of CBS Model}
\label{sec:CBSmodel}
The total cost is a function of $\pi_{{N-1}_A}$, $P_B$ (or $P_A$), and $L_Q$. We obtain these quantities by utilizing our decomposition method. The first task is to decompose the full MC into five (partially overlapping) subchains: $A_1=\{0_{A},1_{A} ,\cdots,n-1_{A}\}$, $A_2=\{n_{A},n+1_{A} ,\cdots,N-1_{A}\}$, $A_3=\{N-1_A, n+1_B\}$, $A_4=\{n+1_{B},n+2_{B} ,\cdots,N_{B}\}$, $A_5=\{N+1_{B},N+2_{B},\cdots\}$. We set $J^+=\{1,2,3,4,5\}$ and $J^-=\{3\}$ to satisfy Equation (\ref{Condition for a set J}). (Note: Subchain 3 can be dropped from both $J^+$ and $J^-$. However, we include subchain 3 since it is convenient to utilize subchain 3 in the analysis.) Let the reference state of the full MC be state $N-1_{A}$. Let $N_Q$ be the operator for the number of vehicles in a queue. By plugging $f(X)=1$, $I_B(\doteq I(k \in  A_4 \cup A_5))$, and $N_Q$ into Theorem \ref{total expectation}, three quantities of interest are obtained by summing up indicators of subchains:

(i) $\pi _{N-1_{A}}$:
\begin{equation*}
\frac{1}{\pi _{N-1_{A}}}=\frac{1}{\beta_{N-1_{A},n-1_{A}} \cdot \pi _{{n-1_{A}}}^1}+\frac{1}{\pi _{{N-1_{A}}}^2}+\frac{1}{\beta_{N-1_{A},n+1_{B}} \cdot \pi _{{n+1_{B}}}^4}+\frac{1}{\beta_{N-1_{A},N+1_{B}} \cdot \pi _{{N+1_{B}}}^5}.
\end{equation*}

(ii) $P_B$:
\begin{equation*}
\frac{P_B}{\pi _{N-1_{A}}}=\frac{E[I_B]}{\pi _{N-1_{A}}}=\frac{1}{\beta_{N-1_{A},n+1_{B}} \cdot \pi _{{n+1_{B}}}^4}+\frac{1}{\beta_{N-1_{A},N+1_{B}} \cdot \pi _{{N+1_{B}}}^5}.
\end{equation*}

(iii) $L_Q$:
\begin{equation*}
\frac{L_Q}{\pi _{N-1_{A}}}=\frac{E[N_Q]}{\pi _{N-1_{A}}}=\frac{L_Q^1}{\beta_{N-1_{A},n-1_{A}} \cdot \pi _{{n-1_{A}}}^1}+\frac{L_Q^2}{\pi _{{N-1_{A}}}^2}+\frac{L_Q^4}{\beta_{N-1_{A},n+1_{B}} \cdot \pi _{{n+1_{B}}}^4}+\frac{L_Q^5}{\beta_{N-1_{A},N+1_{B}} \cdot \pi _{{N+1_{B}}}^5}.
\end{equation*}

Our next task is to analyze the five decomposed subchains. Derivation of closed-form solutions for the subchains is straightforward since these subchains can be analyzed independently by truncation (subchains $A_1$ and $A_5$) or termination (subchains $A_2$, $A_3$, and $A_4$) (see Appendix \ref{sec: EC-CBS} for the derivation in detail). For notational convenience, let $X$ be a Poisson random variable with parameter $\lambda /\mu :E[X]=var(X)=\lambda /\mu $; $\Pr \{ X = s\}$ and $\Pr \{ X\le s\}$ represent Poisson probability mass function and cumulative distribution function, respectively. We denote $s=c-e$, $\rho =\lambda/(s\mu)$, $\omega=1/\rho=s\mu/\lambda$, and $\eta=\lambda/(c\mu)$. We assume that $n, N, c, e$, and $s$ are all integers that satisfy $N>n \geq c$, $e>0$, and $s=c-e>0$. Define four functions as follows (Appendices \ref{sec: EC-mmsk} and \ref{sec: EC-mm1r}):
$$f^1(k,\omega)\doteq \frac{\Pr \{ X\le s\} \omega^{k-s} }{\Pr \{ X=s\} } +\frac{1-\omega^{k-s}}{1-\omega},g^1(k,\omega)\doteq \frac{(k-s)-(k-s+1)\omega +\omega ^{k-s+1} }{(1-\omega)^{2}},$$
$$f^2(k,\omega)\doteq \frac{k+1}{1-\omega}-\frac{\omega(1-\omega^{k+1})}{(1-\omega)^2}, \text{ and } g^2(k,\omega)\doteq \frac{k(k+1)}{2(1-\omega)} - \frac{(k-(k+1)\omega+\omega^{k+1})\omega}{(1-\omega)^3}.$$

Performance indicators of subchains are summarized as follows:

(i) subchain $A_1$:
$$\frac{1}{\pi _{n-1_{A} }^1} =f^1(n-1,\omega) \text{ and } \frac{L_Q^1}{\pi _{n-1_{A}}^1} =g^1(n-1,\omega).$$

(ii) subchain $A_2$:
$$\frac{1}{\pi _{N-1_{A}}^2} =f^2(N-n-1,\omega) \text{ and } \frac{L_Q^2}{\pi _{N-1_{A}}^2} =(n-s)f^2(N-n-1,\omega)+g^2(N-n-1,\omega).$$

(iii) subchain $A_4$:
$$\frac{1}{\pi_{n+1_{B}}^4} =f^2(N-n-1,\eta) \text{ and } \frac{L_Q^4}{\pi_{n+1_{B}}^4} =(N-c)f^2(N-n-1,\eta)-g^2(N-n-1,\eta).$$

(iv) subchain $A_5$:
$$\frac{1}{\pi _{N+1_{B}}^5} =\frac{1}{1-\eta} \text{ and } \frac{L_Q^5}{\pi_{N+1_{B}}^5} =\frac{N-c+1}{1-\eta} +\frac{\eta}{\left(1-\eta\right)^{2} }.$$

We can also derive coefficients:
$$\beta_{N-1_{A} ,n-1_{A}} =\frac{1-\omega}{\omega(1-\omega^{N-n} )}, \beta_{N-1_{A},n+1_{B} }=\frac{1}{\eta}, \text{ and } \beta_{N-1_{A} ,N+1_{B}} =\frac{1-\eta}{\eta^{2} (1-\eta^{N-n})}.$$

Combining the above, we obtain the analytical representation for all necessary indicators to calculate the total cost  $C_{staff}+C_{wait}+C_{switch}$:

(i) $\pi _{N-1_{A}}$:
\begin{equation*}
\frac{1}{\pi _{N-1_A}} =\left(\frac{1}{1-\omega} +\frac{\eta}{1-\eta} \right)(N-n) -\frac{\omega^{n-s} (1-\omega^{N-n} )}{1-\omega} \left(\frac{1}{1-\omega} -\frac{\Pr \{ X\le s\} }{\Pr \{ X=s\} } \right).
\end{equation*}

(ii) $P_B$:
\begin{equation*}
\frac{P_B}{\pi _{N-1_{A}}}=\frac{\eta}{1-\eta} (N-n).
\end{equation*}

(iii) $L_Q$:
\begin{equation*}
\frac{L_Q}{\pi _{N-1_{A}}}=\left(\frac{1}{1-\omega} +\frac{\eta}{1-\eta} \right)(N-n) \left(\frac{N+n+1}{2} - \frac{\lambda}{\mu}-\frac{1}{1-\omega} +\frac{\eta}{1-\eta} \right).
\end{equation*}

\begin{figure}[h]
\FIGURE
{\includegraphics*[scale=0.5]{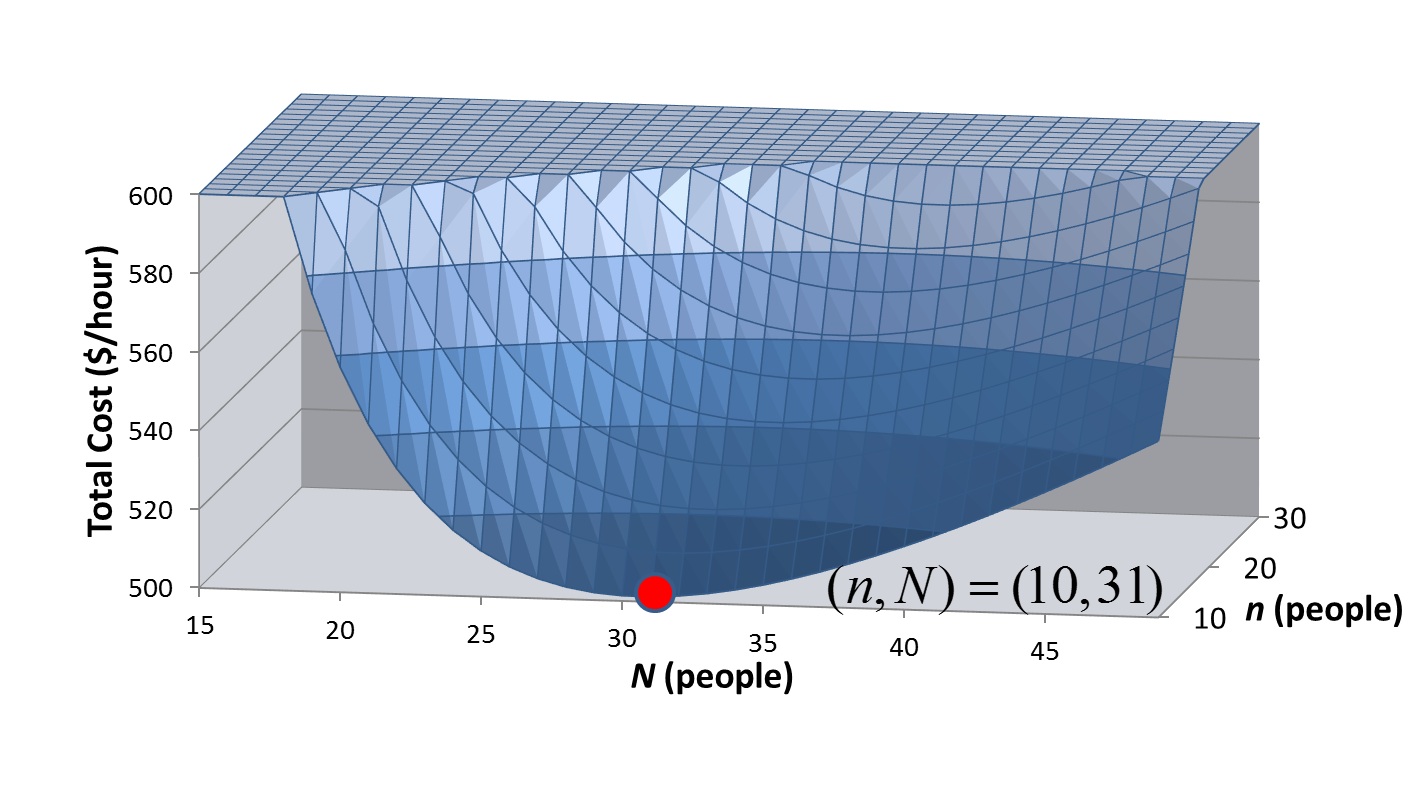}}
{Total cost of operation as a function of $(n,N)$ thresholds with $\lambda=600\text{/hour}, \mu=65\text{/hour}, c=10, e=s=5, c_g=\$20\text{/hour}, c_w=\$10\text{/hour}, \text{ and } c_s=\$50\text{/cycle}$.\label{fig:TotalCost-CBS}}
{}
\end{figure}

Using these expressions, we obtain a closed-form representation of the total cost, which is easy to evaluate in an Excel spreadsheet. Figure~\ref{fig:TotalCost-CBS} is the \emph{exact} cost surface for all possible $(n,N)$ combinations, from which we obtain the optimal CBS policy with thresholds $(n,N)=(10,31)$ for the parameters given in the caption of Figure~\ref{fig:TotalCost-CBS}. This example illustrates the simplicity of our approach: We can derive analytical representations by a simple summation of performance indicators of subchains. For example, if we want to include customer abandonment in subchain 5 in this CBS model, we only need to replace the terms originated from subchain 5 with the new terms in the summation. Or if we consider adding a third threshold (or more), we only need to add terms corresponding to the subchains representing new thresholds in the summation: The complexity of the computation following our method is increased only linearly with respect to the number of subchains added.
\textcolor{black}{\subsection{Analysis of M\textsubscript{t}/M\textsubscript{t}/1 Queue}
\label{sec:mtmt1queue}
We analyze a single-server Markov-modulated MC, an M\textsubscript{t}/M\textsubscript{t}/1 queue, whose arrival and service rates are time-dependent (see Figure \ref{fig:mtmt1queue}).
\begin{figure}[h]
\FIGURE
{\includegraphics*[scale=0.4]{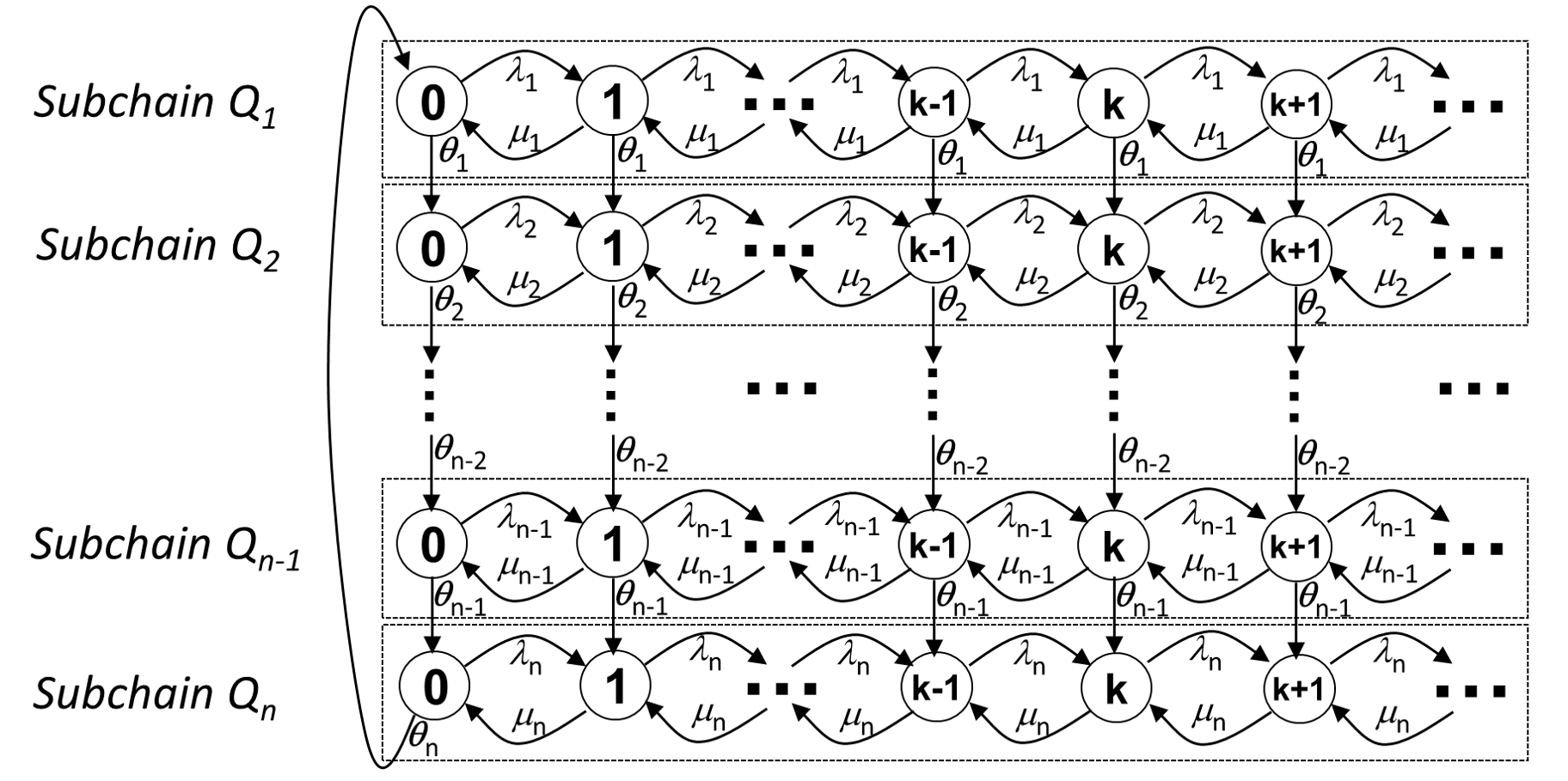}}
{M\textsubscript{t}/M\textsubscript{t}/1 Queue.\label{fig:mtmt1queue}}
{}
\end{figure}
We show the procedure to derive the z-transform (moment-generating function) of the number of customers (jobs) in this system. Let the states of the M\textsubscript{t}/M\textsubscript{t}/1 queue be composed of $Q_1=\{0_1, 1_1, 2_1, \cdots\}, Q_2=\{0_2, 1_2, 2_2, \cdots\},\cdots, Q_n=\{0_n, 1_n, 2_n, \cdots\}$, where a state $\ell_k$ represents $\ell$ customers in the $Q_k$ queue (or a state at \emph{phase} $k$ and \emph{level} $\ell$) in the system. We assume that the system transitions from phase $k$ to $k+1$ (i.e., $Q_k$ to $Q_{k+1}$) with the inter-transition rate $\theta_k$ for $k \in \{1,\cdots,n-1\}$, and then from phase $n$ to zero ($Q_n$ to $Q_0$), which occurs only from states $0_n$ to $0_1$ with rate $\theta_n$. This M\textsubscript{t}/M\textsubscript{t}/1 queue is often used to model an operational problem involving machine deterioration and replacement/repair, where the replacement corresponds to the last transition from $Q_n$ to $Q_0$, which occurs only when there are no jobs to process. As is typical, we assume our single machine maintains the same processing speed for exponential time when the system is in the same phase (i.e., subchain $Q_k$). The machine decreases its processing speed as the system moves to the next phase. When the system is in the last phase ($Q_n$) and a job is cleared (state $0_n$), the machine is renewed with a Poisson rate.}

\textcolor{black}{Let the random variable $X$ represent the number of customers in the system. Our goal is to find the z-transform $P(z)=E[z^X]$. We present termination schemes that conserve the correct distributions of decomposed subchains, reveal insights into the dependencies among decomposed subchains, and explain the procedure to derive the z-transform. Our focus in this example is to understand the analytical properties and interrelations within a complex MC composed of multiple simple subchains, which are often disregarded when solving it numerically.}
\textcolor{black}{\subsubsection{Decomposition Procedure: Analysis of Two Interconnected Queues}
\label{sec:stackedqueue}
We first consider a subproblem of analyzing a part of the M\textsubscript{t}/M\textsubscript{t}/1 queue: two stacked subchains $Q_\pi$ and $Q_p$ (Figure \ref{fig:stackedqueue}).
\begin{figure}[h]
\FIGURE
{\includegraphics*[scale=0.4]{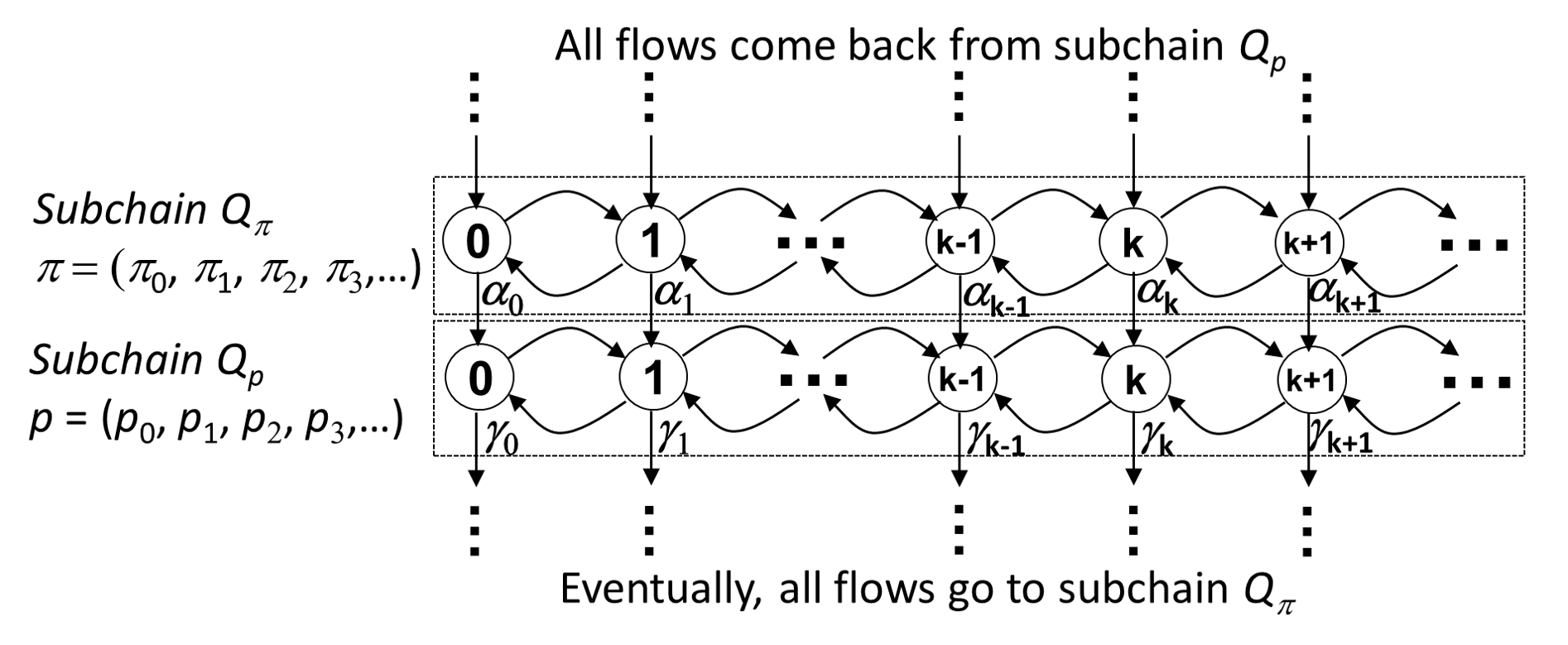}}
{Two Stacked Queues.\label{fig:stackedqueue}}
{}
\end{figure}
In this model, we assume that there is a direct transition from each state in $Q_\pi$ to the corresponding state in $Q_p$. We also assume that all transitions from states in $Q_p$ reach a communication class that includes $Q_\pi$ (thus, transitions from states in $Q_p$ to states in $Q_{\pi}$ are not necessarily direct.) We denote the states in $Q_p$ as $\{0_p,1_p,\cdots,k_p,\cdots\}$ and the states in $Q_\pi$ as $\{0_\pi,1_\pi,\cdots,k_\pi,\cdots\}$, and the corresponding distributions as $p=(p_0,p_1,\cdots,p_k,\cdots)$ and $\pi=(\pi_0,\pi_1,\cdots,\pi_k,\cdots)$, respectively. To analyze $Q_p$, we assume its MC structure (after truncating all inter-transitions between subchains) is an M/M/1 queue with arrival rate $\lambda$ and service rate $\mu$. The structure of $Q_\pi$ is arbitrary except we assume that its distribution is $\pi$. Our goal is to find the z-transform of $p$ given ${\pi}$.}

\textcolor{black}{We introduce two key parameters, $\hat{\pi}_k$ and $\hat{p}_k$, which play an important role in our decomposition method. Using the average rate of inflow (into $Q_p$) $\bar{\alpha}$ and the average rate of outflow (from $Q_p$) $\bar{\gamma}$, we define the parameters as
\begin{equation}
\label{proportion}
\hat{\pi}_k\doteq\alpha_k \pi_k/\bar{\alpha},\,\text{where } \bar{\alpha}=\sum_{i=0}^\infty \alpha_i\pi_i, \quad
\hat{p}_k\doteq\gamma_k p_k /\bar{\gamma},\, \text{where } \bar{\gamma}\doteq\sum_{k=0}^\infty \gamma_i p_i.
\end{equation}
A set of parameters $\hat{\pi}=(\hat{\pi}_0,\hat{\pi}_1,\cdots,\hat{\pi}_k,\cdots)$ represents a proportion of inflow to all states in $Q_p$; similarly, $\hat{p}=(\hat{p}_0,\hat{p}_1,\cdots,\hat{p}_k,\cdots)$ represents a proportion of outflow from all states in $Q_p$. Since $\hat{\pi}$ and $\hat{p}$ are proportions, they are, by definition, normalized (and can be treated just like stationary probabilities): $\sum_{i=0}^\infty \hat{\pi}_i=\sum_{i=0}^\infty \hat{p}_i=1$.}

\textcolor{black}{Following our termination scheme (Corollary \ref{special termination}), we can conserve partial flow at all states in $Q_p$ if all transitions out from $Q_p$ are redirected back into $Q_p$ such that the flow to each state $k_p \in Q_p$ is in proportion to the inflow to $k_p$, $\hat{\pi}_k$. Specifically, we terminate $Q_p$ by adding transitions from state $k_p$ to state $k'_p$ with rate $\Delta^p_{k,k'}=\gamma_k \hat{\pi}_{k'}$ for all $k$ and $k'$ (including self-transitions, although these do not affect the final result); Figure \ref{fig:case-a} shows how we terminate $Q_p$ to maintain its correct distribution after decomposition. This termination scheme may look complicated, but an appropriate termination scheme often takes a much simpler form, depending on the inter-transition rates $\alpha_k$ and $\gamma_k$; we present a termination scheme for each specific case below. Note, though, that for the purpose of deriving, for example, the z-transform, all we need is the general termination scheme indicated in Figure \ref{fig:case-a}; we do not need to know a specific termination scheme for each case.
\begin{figure}
    \caption{Termination Schemes.}
    \label{fig:}
    \centering
       \begin{subfigure}[]{0.48\textwidth}
     \centering
        \includegraphics[scale=0.32]{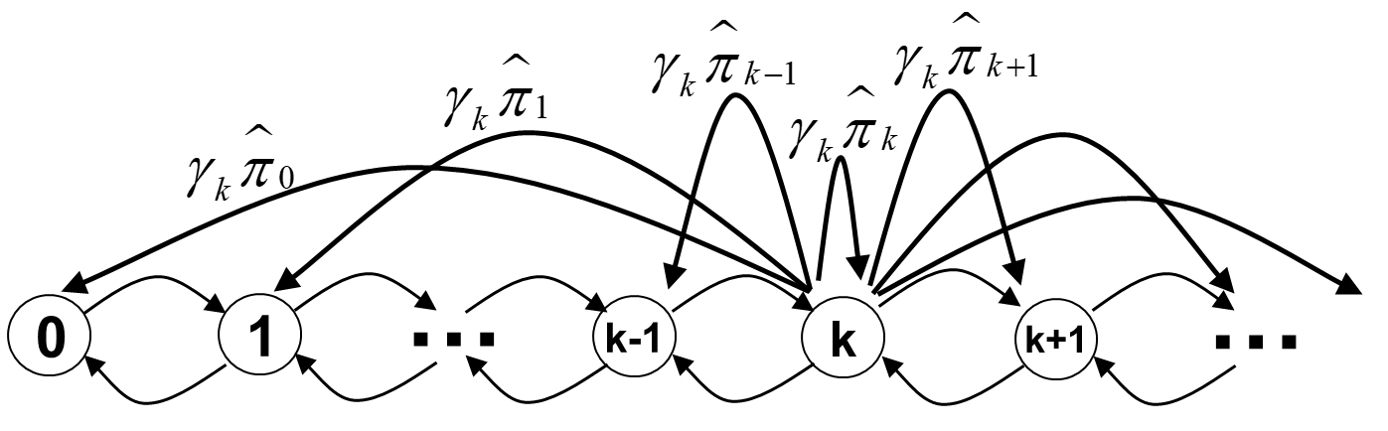}
        \caption{\footnotesize General Case.}
        \label{fig:case-a}
   \end{subfigure}
   \begin{subfigure}[]{0.48\textwidth}
     \centering
        \includegraphics[scale=0.32]{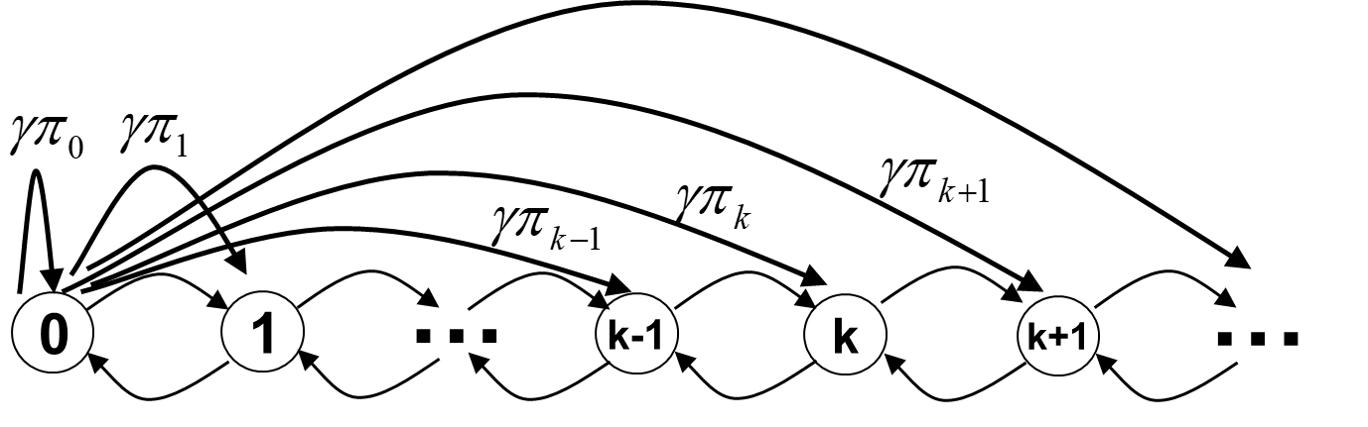}
        \caption{\footnotesize Subchain $Q_n$.}
        \label{fig:case-b}
    \end{subfigure}
    \begin{subfigure}[]{0.48\textwidth}
     \centering
        \includegraphics[scale=0.32]{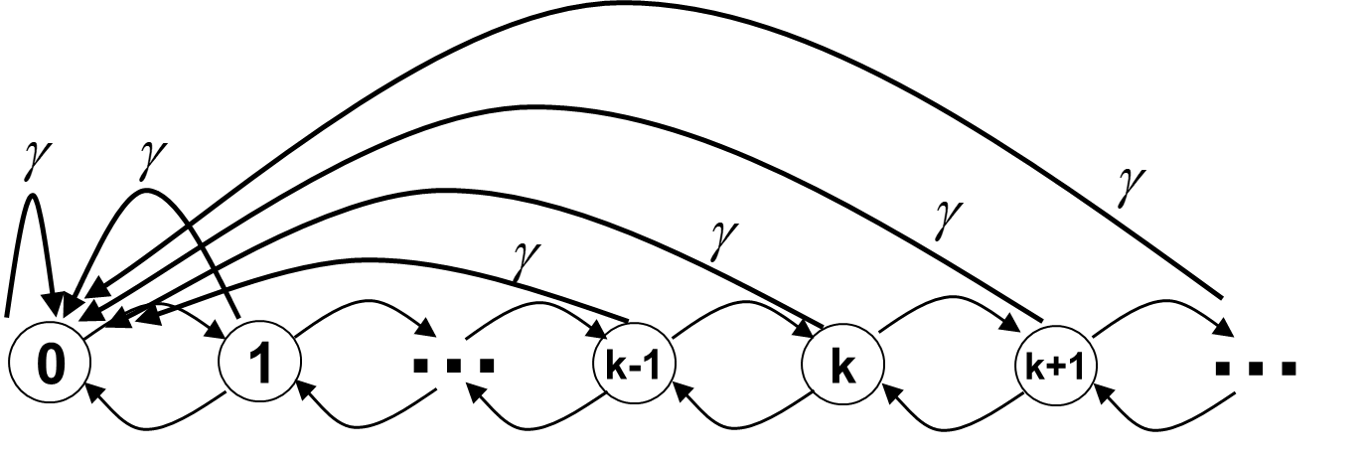}
        \caption{\footnotesize Subchain $Q_1$.}
        \label{fig:case-c}
   \end{subfigure}
   \begin{subfigure}[]{0.48\textwidth}
     \centering
        \includegraphics[scale=0.32]{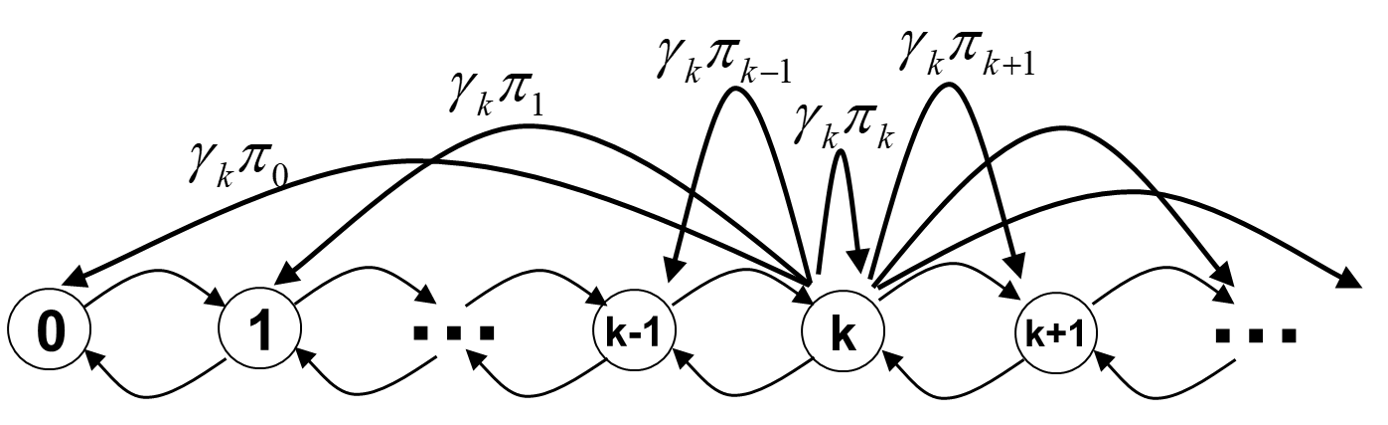}
        \caption{\footnotesize Subchains $Q_2, \cdots, Q_{n-1}$.}
        \label{fig:case-d}
    \end{subfigure}
\begin{flushleft}
\hspace{1.5cm} \footnotesize \emph{Note}: In subfigures (a) and (d), only transitions out from state $k$ are indicated.
\end{flushleft}
\end{figure}}

\textcolor{black}{The equilibrium equations for $Q_p$ are, including self-transitions,
\begin{equation*}
\lambda p_0+\gamma_0 p_0 = \left(\sum_{i=0}^\infty \gamma_i \pi_i\right) \hat{\pi}_0 + \mu p_1, \quad   k=0,
\end{equation*}
\begin{equation*}
(\lambda+\mu) p_k+\gamma_k p_k = \left(\sum_{i=0}^\infty \gamma_i \pi_i\right) \hat{\pi}_k + \mu p_{k+1}+\lambda p_{k-1}, \quad   k\ge 1.
\end{equation*}}

\textcolor{black}{Define z-transforms as follows: $P(z)\doteq\sum_{k=0}^\infty p_k z^k$, $\hat{P}(z)\doteq\sum_{k=0}^\infty \hat{p}_k z^k$, ${\Pi}(z)\doteq\sum_{k=0}^\infty {\pi}_k z^k$, and $\hat{\Pi}(z)\doteq\sum_{k=0}^\infty \hat{\pi}_k z^k$. (Note: All of these take a value of 1 when z=1 because distributions and proportions are normalized.) Also, recall that $\gamma_k p_k = \bar{\gamma} \hat{p}_k, \forall k$, and $\bar{\gamma}=\sum_{i=0}^\infty \gamma_i p_i$  from Equation \eqref{proportion}. To utilize z-transforms, we multiply $z^k$ to the $k\textsuperscript{th}$ ($k\ge1$) equation and sum them up:
\begin{equation*}
(\lambda+\mu) \sum_{k=1}^\infty p_k z^k+\bar{\gamma} \sum_{k=1}^\infty \hat{p}_k z^k = \bar{\gamma} \sum_{k=1}^\infty \hat{\pi}_k z^k + \frac{\mu}{z} \sum_{k=1}^\infty p_{k+1}z^{k+1}+\lambda z \sum_{k=1}^\infty  p_{k-1} z^{k-1},
\end{equation*}
which can be rewritten as
\begin{equation*}
(\lambda+\mu) [P(z)-p_0]+\bar{\gamma} [\hat{P}(z)-\hat{p}_0] = \bar{\gamma}[\hat{\Pi}(z) - \hat{\pi}_0] + \frac{\mu}{z} [P(z)-p_0-p_1 z]+\lambda z P(z).
\end{equation*}
Note that the 0\textsuperscript{th} equation can be rewritten as
\begin{equation*}
\lambda p_0 + \bar{\gamma} \hat{p}_0 = \bar{\gamma} \hat{\pi}_0 + \mu p_1.
\end{equation*}
Thus, we obtain
\begin{equation}
\label{MMPP}
P(z)=\frac{(1-z)p_0+\frac{\bar{\gamma}}{\mu} z\left[\hat{P}(z)-\hat{\Pi}(z)\right]}{(1-\frac{\lambda}{\mu}z)(1-z)},
\end{equation}
which holds for a general class of two stacked queues when one is an arbitrary queue with distribution $\pi$ and the other is an M/M/1 queue with distribution $p$, with inter-transition rates $\alpha_i$ and $\gamma_i$. Note that Equation \eqref{MMPP} holds even if inter-transitions are generalized (e.g., $\alpha_i$ is extended to $\alpha_{i,j}, i \in Q_\pi, j \in Q_p$), in which case we just have to re-define proportions $\hat{\pi}$ and $\hat{p}$ appropriately.}

\textcolor{black}{Equation \eqref{MMPP} indicates that $Q_\pi$ impacts $Q_p$ only via $\hat{\Pi}(z)$; thus, $Q_p$ is independent from $Q_\pi$ only when $\hat{\Pi}(z)$ does not depend on $\pi$. This occurs when a single $\alpha_k=\alpha>0$ and other $\alpha_{k'}=0$ for $k' \neq k$; in this case we obtain $\hat{\Pi}(z)=z^k$, which is independent from $\pi$. In particular, if $k=0$ (i.e., an inflow to $Q_p$ is only observed at state $0_p$), then $\hat{\Pi}(z)=1$, which is obviously independent from $\pi$. (We call this special case \emph{single-channel} and study it further below.) For any other cases, $Q_p$ depends on $Q_\pi$.}

\textcolor{black}{From the above analysis, if an MC is composed of a series of stacked queues as shown in Figure \ref{fig:stackedqueue}, such an MC can be solved following one of the three possible scenarios (\emph{trichotomy}): i) all subchains can be solved independently, ii) one subchain can be solved independently and others are obtained subsequently, and iii) all/some need to be solved simultaneously. Many interesting queueing models belong to the first two cases due to their special structures, and thus yield analytical solutions. The M\textsubscript{t}/M\textsubscript{t}/1 queue we study is one of them; based on the inter-transition rates of the model, it is evident that this queue can be solved sequentially (the second case). Specifically, the first subchain $Q_1$ is solved independently, $Q_2$ is derived given the result of $Q_1$, $Q_3$ is derived given the result of $Q_2$, and we repeat this process down to the last subchain $Q_n$.}

\textcolor{black}{To derive the detailed relationships among subchains, observe that the z-transforms $\hat{\Pi}(z)$ and $\hat{P}(z)$ play a major role in Equation \eqref{MMPP}. These z-transforms may be simplified depending on the inter-transition rates $\alpha_i$ and $\gamma_i$, as we see in the discussion above. We only consider two special cases that are related to the M\textsubscript{t}/M\textsubscript{t}/1 model we study; there are other special cases that can lead to simpler results (e.g., a single-channel case where an inflow is into a single state $k_p, k\ge 1$).}

\textcolor{black}{(a) Single-Channel (SC) Case: If an inflow to $Q_p$ is into a single state $0_p$, then $\alpha_0=\alpha>0, \alpha_i=0, \forall i \ge 1$, and thus $\hat{\Pi}(z)=1$. If an outflow from $Q_p$ is from a single state $0_p$, then $\gamma_0=\gamma>0, \gamma_i=0, \forall i \ge 1$, and thus $\hat{P}(z)=1$.}

\textcolor{black}{(b) Multi-Channel (MC) Case: If an inflow to $Q_p$ is into all states in $Q_p$ with the same transition rate, then $\alpha_i=\alpha,\forall i$, $\bar{\alpha}=\alpha\sum_{i=0}^\infty \pi_i = \alpha$, $\hat{\pi}_k=\alpha \pi_k/\bar{\alpha}=\pi_k, \forall k$, and thus $\hat{\Pi}(z)={\Pi}(z)$. Likewise, if an outflow from $Q_p$ is from all states with the same transition rate, then $\gamma_i=\gamma>0, \forall i$, and thus $\hat{P}(z)={P}(z)$.}

\textcolor{black}{We have four combinations of these special cases, three of which apply to different phases of our M\textsubscript{t}/M\textsubscript{t}/1 system. By plugging appropriate forms of $\hat{\Pi}(z)$ and $\hat{P}(z)$ into Equation \eqref{MMPP}, we obtain the necessary termination scheme as well as the z-transform for each case. The results are summarized in Table \ref{table-MtMt1}, from which we can first derive $P_1(z)$ independently, and then derive $P_2(z)$ to $P_n(z)$ sequentially. The parameter $p_0$ for each case can be derived in a standard manner (see Appendix \ref{sec: EC-MtMt1}).}

\textcolor{black}{Table \ref{table-MtMt1} shows the benefits of our decomposition method. First, it shows the generality of our method; once we obtain the general solution (e.g., Equation \eqref{MMPP}) for a class of MCs, we can obtain the results of special cases in the same class straightforwardly. Second, by explicitly revealing the termination scheme for each case, we can identify an MC that has been analyzed before and can reuse previously obtained results if available. In this paper we solve all cases to show our method, but this is not necessary in practice; for example, when $\hat{\Pi}(z)=\hat{P}(z)=1$, we obtain $\Delta^p_{0,0}=\gamma$, which implies that the necessary termination scheme is a simple truncation and thus $P(z)$ should match the known result for an M/M/1 queue. As another example, when $\hat{\Pi}(z)=1$ and $\hat{P}(z)=P(z)$ hold, we obtain $\Delta^p_{k,0}=\gamma, \forall k$, which implies that the MC is reduced to \emph{A Processor Model with Failures} \citep{nelson1995probability} and thus $P(z)$ should match that for their model. In fact, two of the four special cases in Table \ref{table-MtMt1} are known, and the remaining two can be reused in the future when analyzing other MCs with such subchains. We expect that many complicated MCs are made up of simpler truncated or terminated MCs. Our method can enable a direct reuse of the previously known results if available.}
\begin{table}
\caption{Termination schemes and z-transforms.}
\label{table-MtMt1}
\centering
\begin{tabular}{l|l|l|l}
\multicolumn{1}{l}{}& \multicolumn{1}{c}{Inflow (SC): $\hat{\Pi}(z)=1$}&\multicolumn{1}{c}{Inflow (MC): $\hat{\Pi}(z)=\Pi(z)$}&\\ 
\cline{2-3}
\begin{tabular}[c]{@{}l@{}}Outflow (SC):\\$\hat{P}(z)=1$\end{tabular}       & \begin{tabular}[c]{@{}l@{}}M/M/1 queue\\ $\Delta^p_{0,0}=\gamma$ (Truncation)\\$P(z)=\frac{p_0}{1-\frac{\lambda}{\mu}z},\; p_0=1-\frac{\lambda}{\mu}$\end{tabular} & \begin{tabular}[c]{@{}l@{}}$Q_n$ of M\textsubscript{t}/M\textsubscript{t}/1 queue\\$\Delta^p_{0,k}=\gamma {\pi}_{k}, \forall k$ (Figure \ref{fig:case-b})\\$P(z)=\frac{(1-z)+\frac{{\gamma}}{\mu} z\left[1-{\Pi}(z)\right]}{(1-\frac{\lambda}{\mu}z)(1-z)}p_0,\; p_0=\frac{1-\frac{\lambda}{\mu}}{1+\frac{\gamma}{\mu}{\Pi}'(1)}$\end{tabular} &   \\ 
\cline{2-3}
\begin{tabular}[c]{@{}l@{}}Outflow (MC):\\$\hat{P}(z)=P(z)$\end{tabular} & \begin{tabular}[c]{@{}l@{}}$Q_1$ of M\textsubscript{t}/M\textsubscript{t}/1 queue\\$\Delta^p_{k,0}=\gamma, \forall k$ (Figure \ref{fig:case-c})\\$P(z)=\frac{(1-z)p_0-\frac{\gamma}{\mu}z}{(1-\frac{\lambda}{\mu}z)(1-z)-\frac{\gamma}{\mu}z},\; p_0=\frac{\gamma }{\mu} \frac{r_1}{1-r_1}$\end{tabular}& \begin{tabular}[c]{@{}l@{}}$Q_2,\cdots,Q_{n-1}$ of M\textsubscript{t}/M\textsubscript{t}/1 queue\\$\Delta^p_{k,k'}=\gamma {\pi}_{k'}, \forall k, k'$ (Figure \ref{fig:case-d})\\$P(z)=\frac{(1-z)p_0-\frac{\gamma}{\mu}z {\Pi}(z)}{(1-\frac{\lambda}{\mu}z)(1-z)-\frac{\gamma}{\mu}z},\; p_0=\frac{\gamma }{\mu} \frac{r_1 {\Pi}(r_1)}{1-r_1}$\end{tabular} &  \\ 
\cline{2-3}
\multicolumn{1}{l}{}& \multicolumn{1}{l}{}& \multicolumn{1}{l}{} &
\end{tabular}
\begin{flushleft}
\hspace{1.5cm} \footnotesize \emph{Notes}: SC and MC represent a single-channel case and a multi-channel case, respectively.\\
\hspace{2.5cm} \footnotesize $r_1=\frac{\lambda+\mu+\gamma - \sqrt{(\lambda+\mu+\gamma)^2-4\lambda \mu}}{2\lambda}$. The detail of the derivation is described in Appendix \ref{sec: EC-MtMt1}.
\end{flushleft}
\end{table}
\textcolor{black}{\subsubsection{Aggregation Procedure: Application of Total Expectation Theorem}
We have shown the derivation of the z-transforms of subchains $Q_1,\cdots, Q_n$. The final step is to find the performance measure $P(z)=E[z^X]$, the z-transform of the number of customers $X$ in the system. Since the z-transform is defined as a form of expectation, we can apply our total expectation theorem, which represents the performance measure of the full system as the sum of the subchains' expectations, with a single parameter that plays a role of a normalizing constant. Let the z-transform of the stationary distribution of $Q_k$ be $P_k(z)=E_k[z^X]$. Let $w_k$ be the probability that a system belongs to $Q_k$, and define the ratio of $w_k$ and $w_{k'}$ as $\beta_{k,k'}=\frac{w_k}{w_{k'}}$ (Definition \ref{ratio}). Applying the total expectation theorem to the M\textsubscript{t}/M\textsubscript{t}/1 queue and setting the reference subchain as $Q_1$, we obtain
\begin{equation*}
\label{MMPP7}
\frac{E[z^X]}{w_1}=\sum_{k=1}^n \frac{E_k[z^X]}{\beta_{1,k}}, \text{ or } \frac{P(z)}{w_1}=\sum_{k=1}^n \frac{P_k(z)}{\beta_{1,k}}.
\vspace*{0.1cm}
\end{equation*}
To find $\beta_{1,k}$, notice that the global flow balance condition holds at each subchain: $w_k \bar{\gamma}^k=\text{constant}, \forall k$, where we define $\bar{\gamma}_k$ as the average outflow rate from $Q_k$ following the definition in Equation \eqref{proportion}. In addition, by inspecting the M\textsubscript{t}/M\textsubscript{t}/1 queue, we know $\bar{\gamma}_k=\theta_k,\forall k \in \{1,2,\cdots,n-1\}$ and $\bar{\gamma}_n=\theta_n P_n(0)$. (Note that $P_n(0)$ represents the stationary probability at state $0_n$ for subchain $Q_n$.) Thus, we obtain $\beta_{1,k}=\frac{{\theta}_k}{{\theta}_1}, \forall k \in \{1,2,\cdots,n-1\}$ and $\beta_{1,n}=\frac{{\theta}_n P_n(0)}{{\theta}_1}$. The final representation of $P(z)$ becomes
\begin{equation*}
\frac{P(z)}{w_1}=\sum_{k=1}^n \frac{P_k(z)}{\beta_{1,k}}=\sum_{k=1}^{n-1} \frac{{\theta}_1 P_k(z)}{{\theta}_k}+\frac{{\theta}_1 P_n(z)}{{\theta}_n P_n(0)},\: \text{where }\frac{1}{w_1}=\sum_{k=1}^{n-1} \frac{{\theta}_1}{{\theta}_k}+\frac{{\theta}_1}{{\theta}_n P_n(0)}.
\vspace*{0.1cm}
\end{equation*}
The parameter $w_1$ is derived by setting $z=1$ (which makes all z-transforms equal to 1). This $w_1$ represents the probability that a system belongs to $Q_1$, but at the same time, plays the role of the normalizing constant for the expectation of the full MC (i.e., $P(z)$ of the M\textsubscript{t}/M\textsubscript{t}/1 queue). Here, we impose the normalization condition ($P(1)=P_k(1)=1,\forall k$) not to derive the distribution of the full MC, but to determine a single unknown parameter $w_1$, as is typically done in transform analysis.}

\section{Conclusions}
\label{sec:concl}
\textcolor{black}{In this paper we present a new decomposition method to derive performance measures (expectations of interest) for continuous time MCs. Our method is general and versatile since it relies on the two fundamental principles: an expectation is a linear operator (leading to the total expectation theorem in MC settings), and steady-state distributions of subchains remain unchanged when net average flows are maintained (leading to the partial flow conservation condition). To illustrate our method, we apply our method to several MCs and derive their properties and performance measures.}

\textcolor{black}{This paper only focuses on Markovian queues, but the two principles we utilized to develop our method are applicable to more general, non-Markovian queues as well. In fact, there are existing studies closely related to our second concept; for example, the \emph{rate balance principle} (RBP) has been utilized to derive steady-state distributions and less well-known results of M/G/1 and G/M/c queues \citep{oz2017rate}, and \emph{Queueing and Markov Chain Decomposition} (QMCD) has been utilized to solve the state-dependent M/G/1 queue \citep{abouee2016state}, the state-dependent M/G/1 queue with orbit \citep{baron2018state}, and the multiclass M/G/1 make-to-stock queue \citep{abouee2012strategies}. The study of these queues are beyond the scope of this paper, but we believe our method can serve as an alternative framework to study queueing systems in a non-Markovian setting; we would like to explore such an extension in subsequent papers.}




\bibliographystyle{ormsv080} 
\bibliography{abandonment-ref} 




\ECSwitch


\ECHead{Electronic Companion: Markov Chain Decomposition Based On Total Expectation Theorem}

\section{Proofs}
\label{sec: EC-proofs}
\proof{Proof of Proposition \ref{conservation of distribution}.}
\textcolor{black}{Proof is straightforward using Definition \ref{ratio} and is omitted.}
\endproof

\proof{Proof of Theorem \ref{total expectation}.}
\textcolor{black}{We prove a variant of Equation \eqref{tetequation}, whose right and left hand sides (RHS and LHS, respectively) are multiplied by $\pi_k$. Using Definition \ref{ratio}, Proposition \ref{conservation of distribution}, and Equation \eqref{Condition for a set J},
\begin{eqnarray*}
 RHS \cdot \pi_k &=& \sum_{j \in J^+} \frac{E[f(X)|A_j]}{\pi_{(k)}^j} \cdot \pi_k - \sum_{j \in J^-} \frac{E[f(X)|A_j]}{\pi_{(k)}^j} \cdot \pi_k \\
&=& \sum_{j \in J^+} E[f(X)|A_j]P_j - \sum_{j \in J^-} E[f(X)|A_j]P_j\\
&=& \sum_{j \in J^+} \sum_{k \in A_j} f(k) p_k^j P_j - \sum_{j \in J^-} \sum_{k \in A_j} f(k) p_k^j P_j \\
&=& \sum_{j \in J^+} \sum_{k \in A_j} f(k) \pi_k - \sum_{j \in J^-} \sum_{k \in A_j} f(k) \pi_k \\
&=& \sum_{j \in J^+} \sum_{k \in S} I(k \in A_j) f(k) \pi_k - \sum_{j \in J^-} \sum_{k \in S} I(k \in A_j) f(k) \pi_k  \\
&=& \sum_{k \in S} ( \sum_{j \in J^+}  I(k \in A_j) - \sum_{j \in J^-}  I(k \in A_j) ) f(k) \pi_k \\
&=& \sum_{k \in S} f(k) \pi_k = E[f(X)]=LHS \cdot \pi_k. \Halmos
\end{eqnarray*}}
\endproof

\proof{Proofs of Propositions \ref{additivity property} and \ref{linear property}.}
\textcolor{black}{We can easily confirm the results using Theorem \ref{total expectation} and thus omit proofs.}
\endproof

\proof{Proof of Proposition \ref{boundary condition}.}
The equivalence of all conditions in this proposition can be proved straightforwardly. Here, we only show the equivalence between conditions $(a)$ in Proposition \ref{conservation of distribution} and $(d)$ in Proposition \ref{boundary condition}.

$(a) \Leftarrow (d)$: Let a steady-state distribution of the full MC be a row vector $\pi  = \left( {{\pi _{{\mathop{\rm int}} ({A_j})}},{\pi _{\partial {A_j}}},{\pi _{{A_j}^c}}} \right)$, which is divided into three row vectors according to where each state lies. Note that $\pi$  is uniquely determined due to the requirement of ergodicity of the full MC. Let a corresponding transition matrix be $P = \left( {{P_{ij}}} \right),i,j \in \left\{ {{\mathop{\rm int}} ({A_j}),\partial {A_j},{A_j}^c} \right\}$. Since $\pi$ satisfies the equation $\pi  = \pi P$, flow balance equations for the set ${\mathop{\rm int}} ({A_j})$ should follow the equation ${\pi _{{\mathop{\rm int}} ({A_j})}} = {\pi _{{\mathop{\rm int}} ({A_j})}}{P_{{\mathop{\rm int}} ({A_j}),{\mathop{\rm int}} ({A_j})}} + {\pi _{\partial {A_j}}}{P_{\partial {A_j},{\mathop{\rm int}} ({A_j})}}$. The distribution ${\pi _{{\mathop{\rm int}} ({A_j})}}$ should be uniquely solvable in terms of ${\pi _{\partial {A_j}}}$. (If not, there exist multiple solutions for $\pi$.) Therefore, ${\pi _{{\mathop{\rm int}} ({A_j})}} = {\pi _{\partial {A_j}}}{P_{\partial {A_j},{\mathop{\rm int}} ({A_j})}}{\left( {I - {P_{{\mathop{\rm int}} ({A_j}),{\mathop{\rm int}} ({A_j})}}} \right)^{ - 1}}.$

Notice that both ${P_{\partial {A_j},{\mathop{\rm int}} ({A_j})}}$ and  ${P_{{\mathop{\rm int}} ({A_j}),{\mathop{\rm int}} ({A_j})}}$ are not altered by termination (nor by truncation). Hence, the same equation holds for both the full MC and a decomposed subchain $j$ with termination. Since $(d)$ says ${\pi _k} \propto \pi _k^j,\forall k \in \partial {A_j}$, we know ${\pi _k} \propto \pi _k^j,\forall k \in {A_j}$.

$(a) \Rightarrow (d)$: This is obvious because $\partial {A_j} \subset {A_j}$.\Halmos
\endproof
\medskip

\proof{Proof of Lemma \ref{partial flow}}
\textcolor{black}{First, notice that the steady-state probabilities at boundary states for the full MC and a terminated subchain satisfy the following equations:\\
1) Steady-state equations at boundary states for the full MC:
\begin{equation}
\label{equations for the full MC}
{\pi _k}\left(\sum\limits_{k' \in {A_j}}q_{k,k'} + \sum\limits_{k' \in A_j^c}{q_{k,k'}}\right) = \sum\limits_{k' \in {A_j}} {{\pi _{k'}}{q_{k',k}}}  + \sum\limits_{k' \in A_j^c} {{\pi _{k'}}{q_{k',k}}} ,\forall k \in \partial {A_j}.
\end{equation} 
2) Steady-state equations at boundary states for a terminated subchain:
\begin{equation}
\label{equations for a terminated subchain}
\pi _k^j\left( {\sum\limits_{k' \in {A_j}} {{q_{k,k'}}}  + \sum\limits_{k' \in \partial {A_j}} {\Delta {q_{k,k'}^j}} } \right) = \sum\limits_{k' \in {A_j}} {\pi _{k'}^j{q_{k',k}}}  + \sum\limits_{k' \in \partial {A_j}} {\pi _{k'}^j\Delta {q_{k',k}^j}} ,\forall k \in \partial {A_j}.
\end{equation} 
We show the equivalence of conditions $(d)$ in Proposition \ref{boundary condition} and $(f)$ in Lemma~\ref{partial flow} using Equations \eqref{equations for the full MC} and \eqref{equations for a terminated subchain}.}

$(d) \Leftarrow (f)$:
If $(f)$ holds, by combining $(f)$ with Equation (\ref{equations for the full MC}) and eliminating terms with summation over a set $A_j^c$, we obtain $${\pi _k}\left( {\sum\limits_{k' \in Aj} {{q_{k,k'}}}  + \sum\limits_{k' \in \partial Aj} {\Delta {q_{k,k'}^j}} } \right) = \sum\limits_{k' \in Aj} {{\pi _{k'}}{q_{k',k}}}  + \sum\limits_{k' \in \partial Aj} {{\pi _{k'}}\Delta {q_{k',k}^j}} ,\forall k \in \partial {A_j},$$ which is equivalent to Equation (\ref{equations for a terminated subchain}) for $\pi_k^j$. Since the steady-state equations at boundary states for both the full MC and a terminated subchain are equivalent, under the ergodicity assumption for a terminated subchain, we should obtain the same solutions at boundary states (up to a normalization constant). Hence, $(d)$ holds.

$(d) \Rightarrow (f)$:
If $(d)$ holds, then a set of probabilities $\pi_k, \forall k \in \partial A_j$, should satisfy not only Equation (\ref{equations for the full MC}) but also Equation (\ref{equations for a terminated subchain}):
$$\pi _k\left( {\sum\limits_{k' \in {A_j}} {{q_{k,k'}}}  + \sum\limits_{k' \in \partial {A_j}} {\Delta {q_{k,k'}^j}} } \right) = \sum\limits_{k' \in {A_j}} {\pi _{k'}{q_{k',k}}}  + \sum\limits_{k' \in \partial {A_j}} {\pi _{k'}\Delta {q_{k',k}^j}} ,\forall k \in \partial {A_j}.$$
Combining this equation with Equation (\ref{equations for the full MC}) for $\pi_k$ and eliminating terms with summation over a set $A_j$, we obtain $(f)$.
\Halmos
\endproof
\medskip

\proof{Proof of Proposition \ref{termination}}
By inspection, the termination scheme satisfying these two flow conditions satisfies Lemma \ref{partial flow}. Therefore, the boundary distribution and hence the stationary distribution of a decomposed subchain must be conserved.\halmos
\endproof
\medskip

\proof{Proof of Corollary \ref{special termination}}
\textcolor{black}{We can directly confirm that Equation (\ref{general termination}) satisfies Proposition \ref{termination} (1). To show Equation (\ref{general termination}) satisfies Proposition \ref{termination} (2), we use the global flow balance condition: $$\sum\limits_{k \in {\partial{A_j}}} \left( {\pi _k}\sum\limits_{k' \in A_j^c} {{q_{k,k'}}}  - \sum\limits_{k' \in A_j^c} {{\pi _{k'}}{q_{k',k}}}  \right) =0.$$}Finally, since the added transitions (termination) connect all boundary states, this termination scheme fulfills the ergodicity requirement of a decomposed subchain.\halmos
\endproof
\medskip

\section{Example: Four-State CTMC (Extension to the model discussed in \S\ref{sec:four-state})}
\label{sec: extension}
\begin{figure}[h]
\FIGURE
{\includegraphics*[scale=0.4]{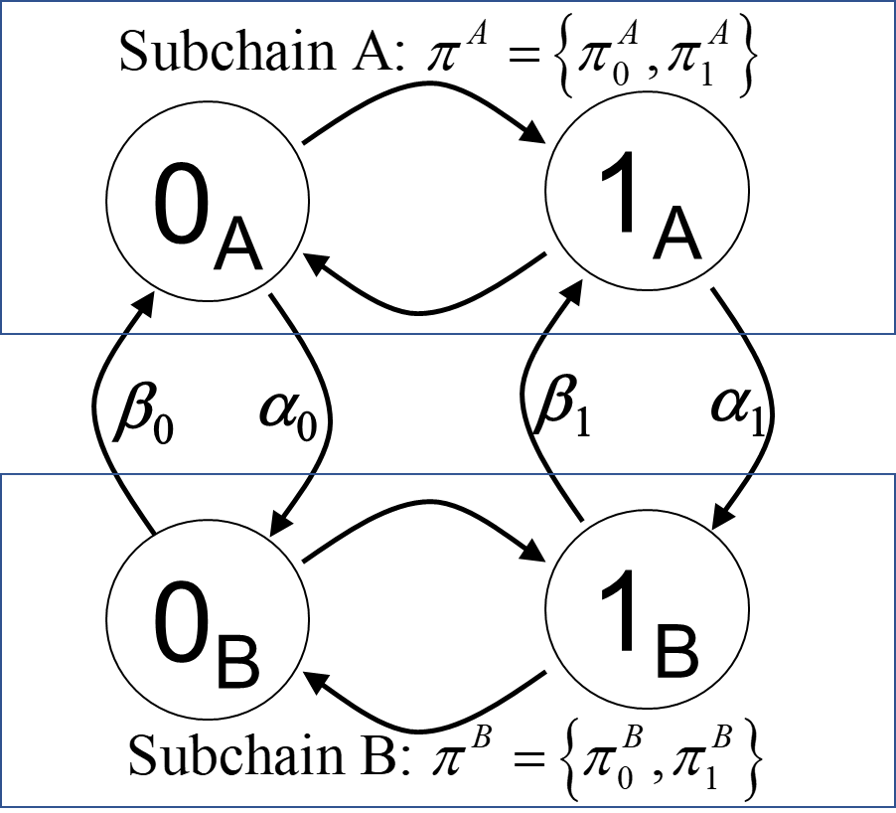}}
{Four-State CTMC.\label{fig:fourstateMC}}
{}
\end{figure}

\textcolor{black}{This example extends the model used in \S\ref{sec:four-state}. We consider again a four-state MC with upper two states $A=\{0_A, 1_A\}$ and lower two states $B=\{0_B, 1_B\}$ (see Figure \ref{fig:fourstateMC}), but this time, we assume that states in $A$ and $B$ are connected via transitions with rate $\alpha_{i}$ ($\beta_{i}$) from a state $i$ in $A$ ($B$) to a corresponding state $i$ in $B$ ($A$, respectively). Let a random variable $X=X_i$, where $i=A$ or $B$. As an example, we derive a performance measure $var(X)$.}
\begin{figure}
    \caption{Trichotomy of Decomposition Analysis.}
    \label{fig:}
    \centering
       \begin{subfigure}[]{0.48\textwidth}
     \centering
        \includegraphics[scale=0.28]{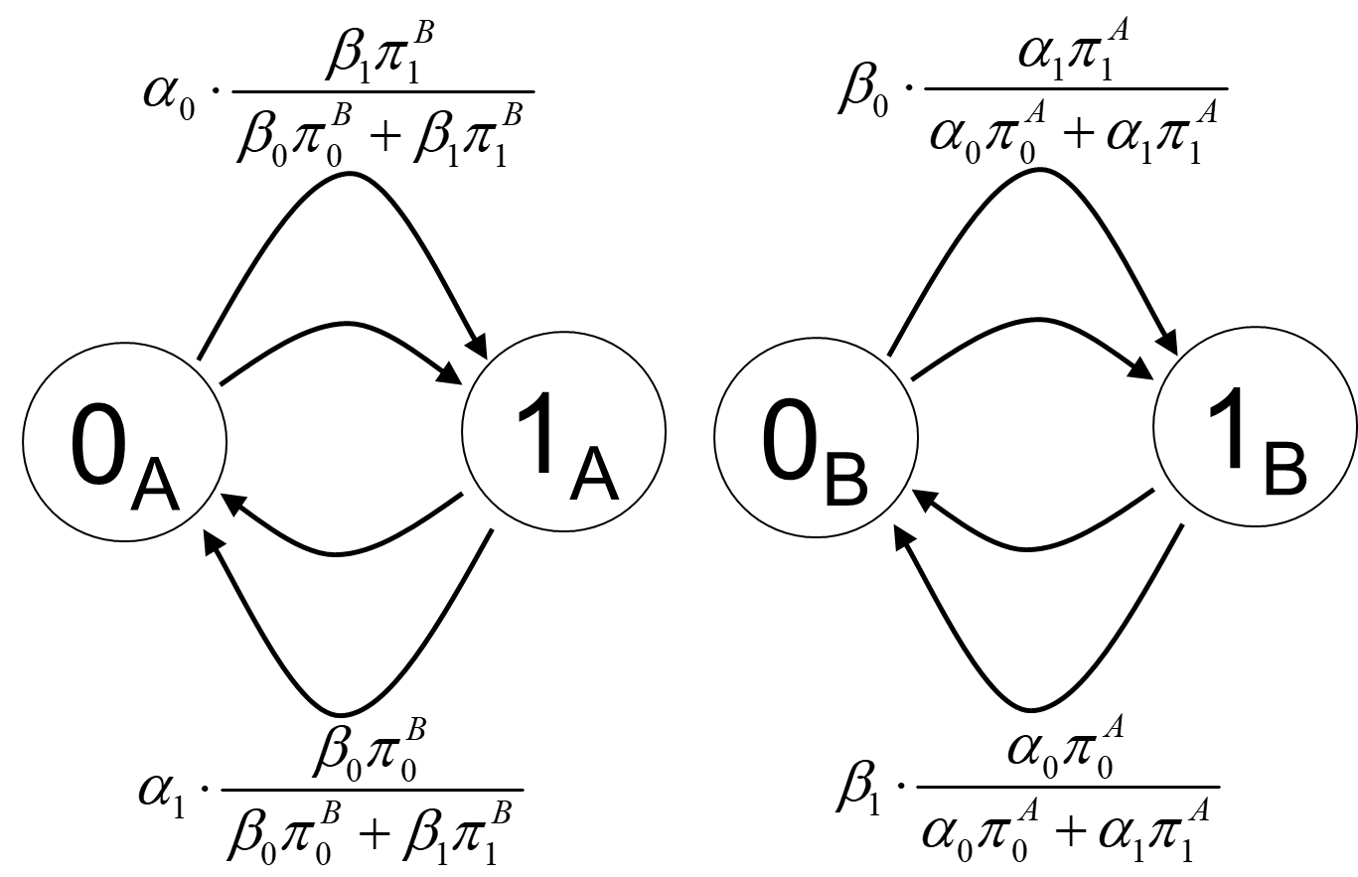}
        \caption{\footnotesize Case I: $\pi^A$ and $\pi^B$ depend on each other.}
        \label{fig:case1}
   \end{subfigure}
   \begin{subfigure}[]{0.48\textwidth}
     \centering
        \includegraphics[scale=0.28]{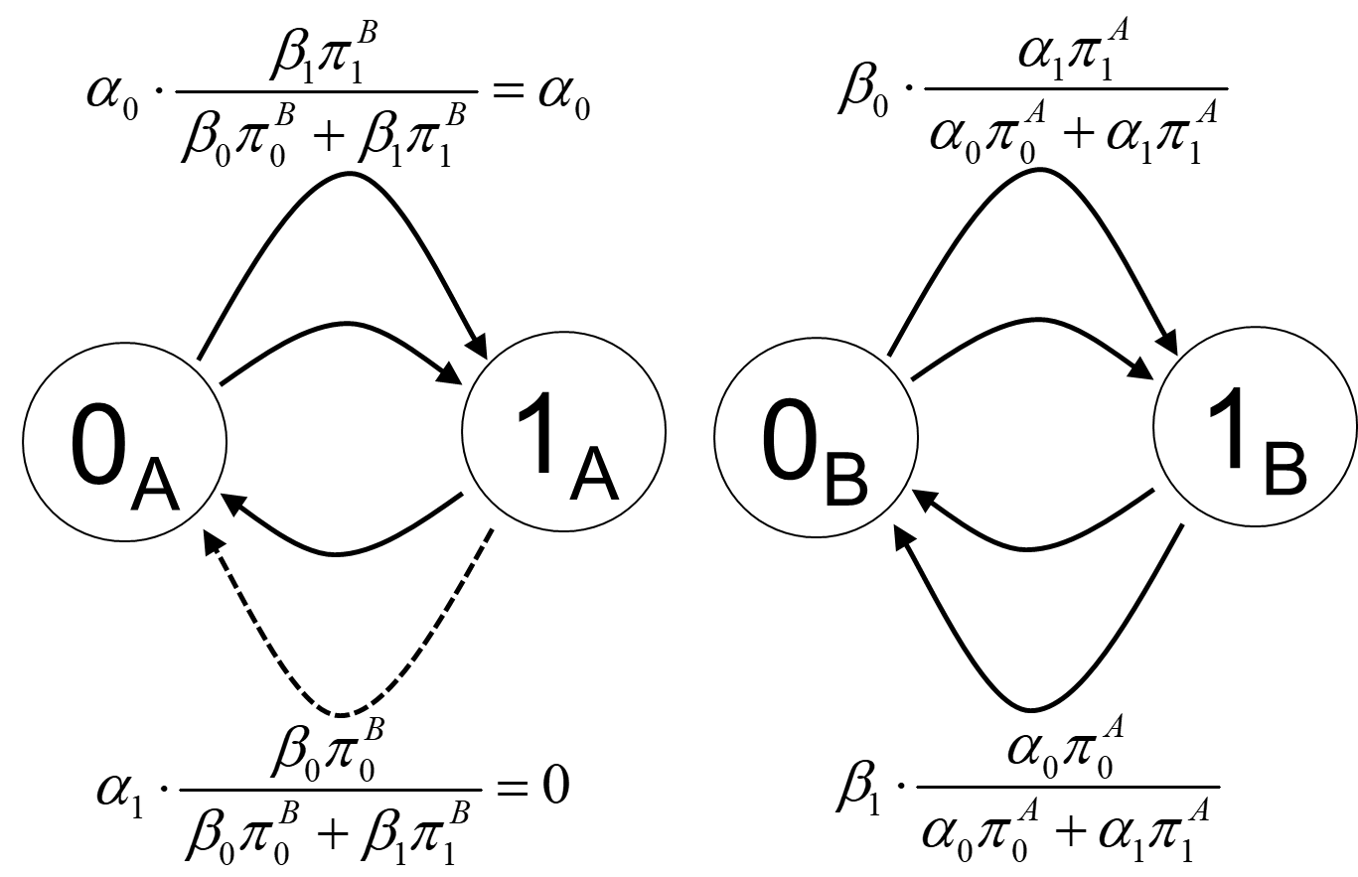}
        \caption{\footnotesize Case II: $\pi^A$ is independent. $\pi^B$ depends on $\pi^A$.}
        \label{fig:case2}
    \end{subfigure}
    \begin{subfigure}[]{0.48\textwidth}
     \centering
        \includegraphics[scale=0.28]{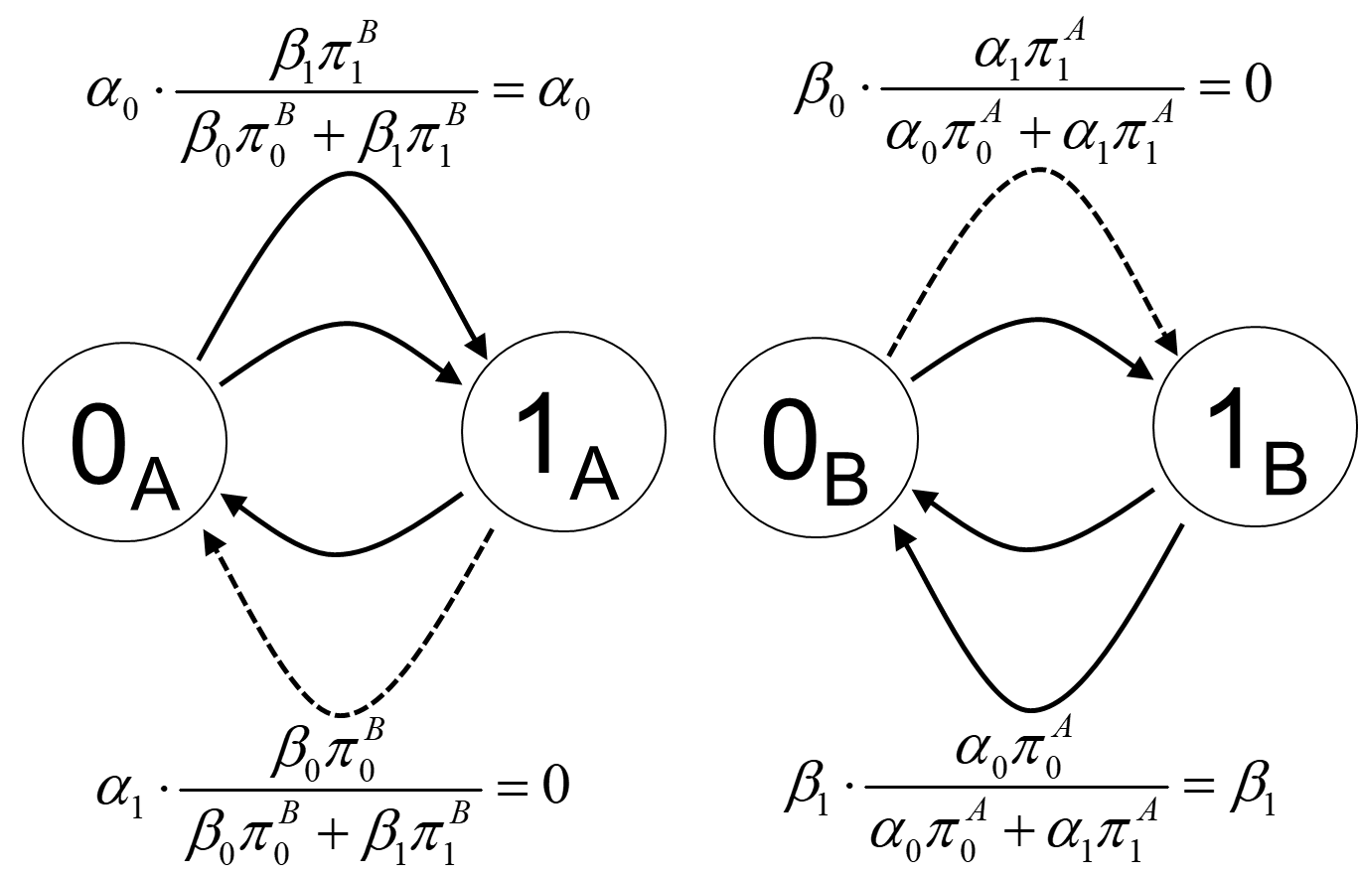}
        \caption{\footnotesize Case IIIa: $\pi^A$ and $\pi^B$ are independent (redirection).}
        \label{fig:case3b}
   \end{subfigure}
   \begin{subfigure}[]{0.48\textwidth}
     \centering
        \includegraphics[scale=0.28]{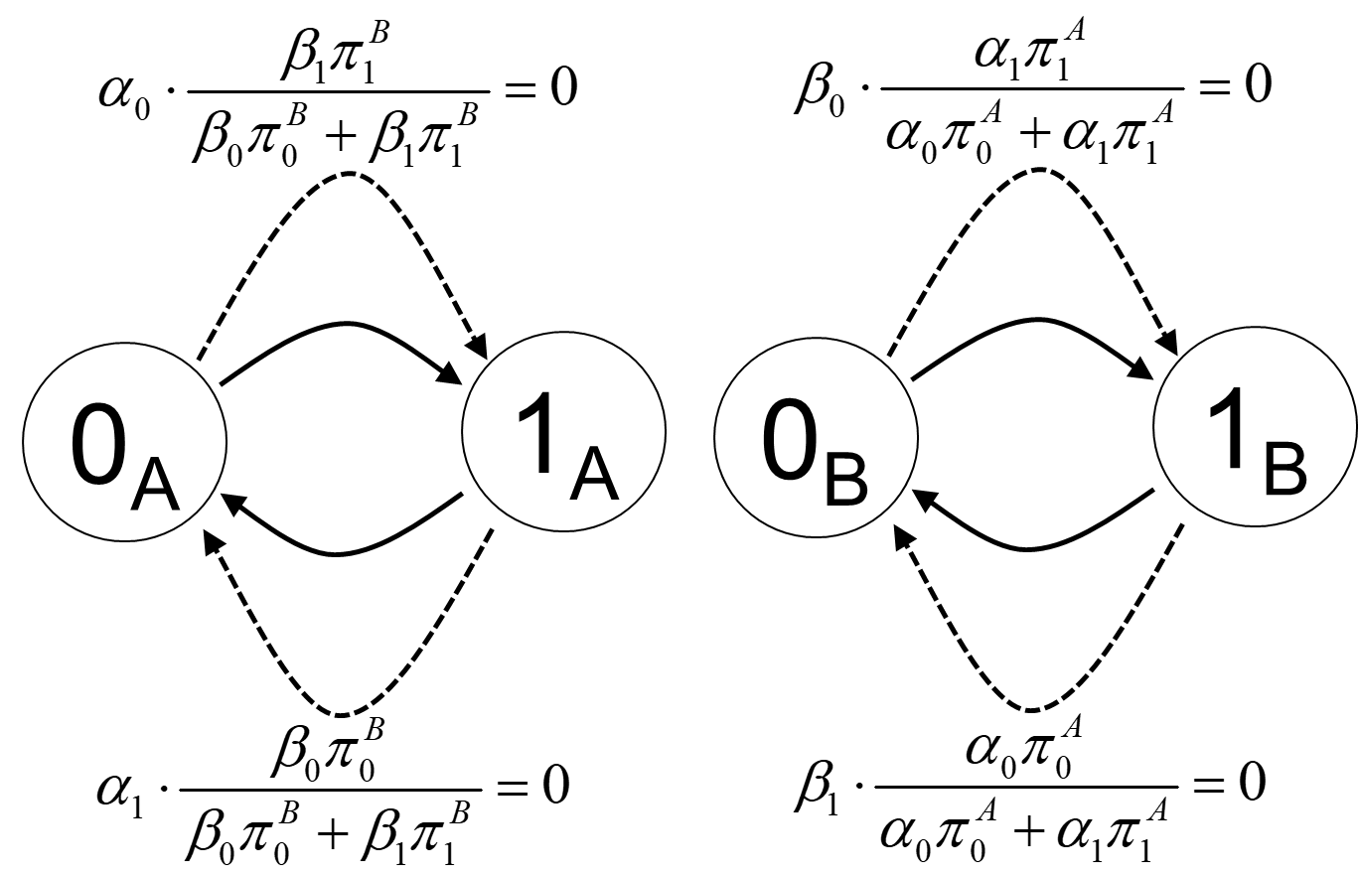}
        \caption{\footnotesize Case IIIb: $\pi^A$ and $\pi^B$ are independent (truncation).}
        \label{fig:case3a}
    \end{subfigure}
\begin{flushleft}
\hspace{.5cm} \footnotesize \emph{Note}: Dashed lines represent vanished transitions.
\end{flushleft}
\end{figure}
\textcolor{black}{Let the distributions of subchains $A$, $B$, and the full MC be $\pi^A=\{\pi^A_0,\pi^A_1\}$, $\pi^B=\{\pi^B_0,\pi^B_1\}$, and $\pi=(\pi_{0A},\pi_{1A},\pi_{0B},\pi_{1B}\}$, respectively. To maintain the correct distributions, i.e., $(\pi^A_0,\pi^A_1) \propto (\pi_{0A},\pi_{1A})$ and $(\pi^B_0,\pi^B_1) \propto (\pi_{0B},\pi_{1B})$, all states of decomposed subchains must conserve partial flow. For this, we \emph{terminate} subchain $A$ by adding transitions $\Delta^A_{i,j}=\alpha_{i} \hat{\pi}^B_j$ from state $i_A$ to state $j_A$, where we define the proportion of inflows (from $B$) to $j_A$ as $\hat{\pi}^B_j=\frac{\beta_{j} \pi^B_j}{\beta_{0} \pi^B_0 + \beta_{1} \pi^B_1}$. This termination scheme \emph{always} works regardless of the structure of subchain $B$ because of the Markov property: the proportion of inflows to $A$ depends only on the distribution of $B$ and not on the transitions made within $B$. This use of the Markov property not only simplifies the termination scheme, but also reveal the impact of $B$ on $A$ explicitly. In our example, a simple observation of termination tells us how to proceed with the analysis. For example, (I) If all $\alpha_i$ and $\beta_i$ are positive, $\pi^A$ and $\pi^B$ need to be solved simultaneously since $\Delta^A_{i,j}$ and $\Delta^B_{i,j}$ depend on $\pi^B$ and $\pi^A$, respectively (Figure \ref{fig:case1}); (II) If $\beta_0=0$ while others are positive, then $\pi^A$ is solved independently and $\pi^B$ is obtained based on $\pi^A$ (Figure \ref{fig:case2}); (IIIa) if $\alpha_1=\beta_0=0$ while others are positive, both $\pi^A$ and $\pi^B$ can be solved for independently by simple redirection of flows (Figure \ref{fig:case3b}); and finally, (IIIb) if $\alpha_0=\beta_0=0$ while others are positive, both $\pi^A$ and $\pi^B$ can be solved for independently by truncation (Figure \ref{fig:case3a}). In general, we observe a \emph{trichotomy} when analyzing MCs using decomposition methods: (I) all/some subchains need to be solved simultaneously; (II) one subchain can be solved independently and others are obtained subsequently; and (III) all subchains can be solved independently by redirection (IIIa), truncation (IIIb), or their combination. Case I typically requires an iterative numerical method or an analysis of simultaneous equations. In contrast, cases II and III can be solved analytically in closed form (with different levels of complexity depending on the model); we show several such examples in this paper (e.g., the M\textsubscript{t}/M\textsubscript{t}/1 queue (Figure \ref{fig:mtmt1queue}) belongs to (II), the CBS model (Figure \ref{fig:MarkovChain-CBS}) belongs to (IIIa), and all birth death queues belong to (IIIb)).}

\textcolor{black}{Now, suppose we obtain $\pi^A$ and $\pi^B$ in the previous step. The next step is to derive a performance measure $var(X)$ utilizing the total expectation theorem in MC settings (Theorem \ref{total expectation}), which in this example takes the form
\begin{equation}
\label{variance}
\frac{E[f(X)]}{w}=E_A[f(X)]+\frac{\alpha}{\beta}E_B[f(X)],
\end{equation}
where we define $\alpha=\sum_k{\alpha_k \pi^A_k}$ ($\beta=\sum_k{\beta_k \pi^B_k}$) as the average transition rate from $A$ to $B$ ($B$ to $A$, respectively), $w$ as a reference set (subchain $A$ in this example) probability, and $f(X)$ as any function of $X$. Equation \eqref{variance} \emph{always} holds regardless of transition rates within each subchain. A normalization condition is part of this relationship, which is utilized to determine $w$. This parameter $w$ is the probability of being in a reference set, but at the same time serves as a normalizing constant in our method.}

\textcolor{black}{Since $var(X)$ is a function of the expected values of $X$ and $X^2$, we first derive $\frac{E[X]}{w}$ and $\frac{E[X^2]}{w}$ by setting $f(X)=X$ and $f(X)=X^2$ in Equation \eqref{variance}, respectively. We can derive the analytical expression: $var(X)=w (E_A[X^2]+\frac{\alpha}{\beta}E_B[X^2])-w^2 (E_A[X]+\frac{\alpha}{\beta}E_B[X])^2$. The parameter $w$ is determined by setting $f(X)=1$ in Equation \eqref{variance}: $\frac{1}{w}=1+\frac{\alpha}{\beta}$. Using subchains' expectations $E_A[X]=E_A[X^2]=\pi^A_1$ and $E_B[X]=E_B[X^2]=\pi^B_1$ and simplifying the result, we obtain $var(X)=p(1-p)$, where $p=\frac{\beta \pi^A_1}{\alpha+\beta}+\frac{\alpha \pi^B_1}{\alpha+\beta}$. This result is anticipated due to the simple example we use here, but illustrates one of the important benefits of our method: we can reveal an interesting relationship between subchains and the original system by taking a divide-and-conquer approach. Following our MC decomposition method, the derivation of such relationships is algorithmic.}

\section{Performance Indicators and Their Relationships}
\label{sec: relationship}
We derive relationships among performance indicators, such as the blocking probability $P_{block}$ (the steady-state probability of the end state of the chain), the queueing probability $P_Q$, and the number of customers in a queue $L_Q$. Such relationships are often difficult to establish, but Theorem \ref{total expectation} enables us to find those relationships without explicitly deriving the performance indicators.

We consider two cases: (1) Subchains share a single state with each other ($|{A_i} \cap {A_j}|=1{\text{ for }}\forall i \ne j$) and (2) Subchains are nested (${A_i} \cap {A_j} = {A_i}{\text{ or }}{A_j}{\text{ for }}\forall i,j$).

\subsection{Subchains Sharing a Single State}
Consider decomposing a general MC into subchains 1 and 2, which share a single state $s$. Define $A_3={A_1} \cap {A_2} = \left\{ s \right\}$. Let $J^+=\{1,2\}$ and $J^-=\{3\}$. As long as proper termination is applied to all subchains, using Theorem \ref{total expectation}, we obtain:
\begin{equation}
\label{general sharing a single state}
\frac{{{E}\left[ {f(X)} \right]}}{{\pi _s}} = \frac{{{E_{{1}}}\left[ {f(X)} \right]}}{{\pi _s^{{1}}}} + \frac{{{E_{{2}}}\left[ {f(X)} \right]}}{{\pi _s^{{2}}}} - f(s), \text{for any } f(X),
\end{equation}
where $\pi_s^3=1$ is used above. This is a general property that always holds regardless of what MC we deal with, how we decompose it, and which function of an MC state we use.

To derive some practically useful properties, consider a simple queueing system which has a chain corresponding to no queue (subchain 1) and a chain corresponding to queue (subchain 2). A state $s$ could be a transitional state shared by both subchains 1 and 2, representing a situation where there is no queue but any new arrival will be put in a queue. (As an example, we can split an M/M/s queue into an M/M/s/s queue and an M/M/1 queue where a state $s$ is shared in the middle.) Denote the queueing probability of this system (i.e., the probability that a new arrival is put in a queue) as $P_Q$, the average number of people in the queue as $L_Q$, and the average number of people in the queue evaluated in subchain 2 as $L_Q^2$. Also, let an operator $N_Q$ represent the number of people in the queue. By applying $f(X)=1$, $I(k \in A_2)$, and $N_Q$ to Equation (\ref{general sharing a single state}), we obtain, correspondingly, 
\begin{equation}
\label{sharing a single state}
\frac{1}{{\pi _s}} = \frac{1}{\pi_s^1} + \frac{1}{{\pi _s^2}} - 1, \frac{P_Q}{\pi_s}=\frac{1}{\pi_s^2}, \text{ and }\frac{L_Q}{\pi_s}=\frac{L_Q^2}{\pi_s^2}.
\end{equation}
Relations among performance indicators are directly obtained from Equation (\ref{sharing a single state}). We denote $\pi_s^1$ as blocking probability $P_{block}$, which represents the probability that a new arrival is blocked to enter a system represented by decomposed subchain 1. By eliminating $\pi_s$, we obtain
\begin{equation}
\label{queueing sharing a single state}
{P_{block}} = \frac{{\pi _s^2{P_Q}}}{{1 - \left( {1 - \pi _s^2} \right){P_Q}}}= \frac{{L_Q}}{{\frac{{L_Q^2}}{\pi_s^2} - \left( {\frac{1}{\pi_s^2} - 1} \right)L_Q}}.
\end{equation}

Equation (\ref{queueing sharing a single state}) holds for any general MC that can be split into no queue and queue subchains, which share a single transitional state $s$. To find the relationship among $P_{block}$, $P_Q$, and $L_Q$, we do not need to derive their actual representations; instead, we need concrete representations of $\pi _s^2$ and $L_Q^2$ in Equation (\ref{queueing sharing a single state}). In particular,

1) If subchain 2 is an M/M/1 queue with $\pi _s^2 = 1 - \rho$ where $\rho<1$, then by using the property $L_Q^2 = \rho / (1 - \rho )$, we obtain
\begin{equation}
\label{P and L relationship-mm1}
{P_{block}} = \frac{{\left( {1 - \rho } \right){P_Q}}}{{1 - \rho {P_Q}}}= \frac{{{{\left( {1 - \rho } \right)}^2}}}{\rho }\frac{{{L_Q}}}{{1 - \left( {1 - \rho } \right){L_Q}}}.
\end{equation}
This relationship is known to hold for Erlang B/C models when subchain 1 is an M/M/s/s queue, and appears in many textbooks (see, for example, \citealt{harchol2013performance}). However, we proved that this relationship holds for \emph{any} subchain 1, not just when subchain 1 is an M/M/s/s queue.
\textcolor{black}{\remark{We observe a similar result in the state-dependent M/G/1 queue when state dependence is for a finite number of states \citep[see Equation (40) in][]{abouee2016state}. They take the \emph{Queueing and Markov Chain Decomposition} (QMCD) approach to decompose the entire system into two subsystems: a state-dependent M/G/1 queue and an auxiliary, state-independent M/G/1 queue (which corresponds to subchain 2 in our analysis). These two subsystems are analyzed separately thanks to the \emph{Level Crossing Theory} (LCT) \citep{brill1977level}, the conservation law of total downcrossing/upcrossing rates for continuous time stochastic processes with continuous state spaces, and then combined to derive the solution, as we did in our analysis. Interested readers may also refer to their analysis on a state-dependent M/G/1/k queue, which is related to our next discussion.}}

2) If subchain 2 is an M/M/1/k queue with $\pi _s^2 = 1 / \sum\limits_{n = 0}^k {{\rho ^n}}$, then by using the property $L_Q^2 = \sum\limits_{n = 0}^k {n{\rho ^n}}\pi_s^2$, we obtain
\begin{equation}
\label{P and L relationship-mmk}
{P_{block}} = \frac{{\left( {1 - \rho } \right){P_Q}}}{{1 - \rho {P_Q} - {\rho ^{k + 1}}\left( {1 - {P_Q}} \right)}}= \frac{{{{\left( {1 - \rho } \right)}^2}}}{\rho }\frac{{{L_Q}}}{{1 - \rho^{k}(1+k(1-\rho)) - \left( {1 - \rho } \right) (1-\rho^k){L_Q}}},
\end{equation}
where we have used the following well-known formulae:
 \[\sum\limits_{n = 0}^k {{\rho ^n}}  = \frac{{1 - {\rho ^{k + 1}}}}{{1 - \rho }} \text{ and } \sum\limits_{n = 0}^k {n{\rho ^n}}  = \rho \frac{\partial }{{\partial \rho }}\left( {\sum\limits_{n = 0}^k {{\rho ^n}} } \right) = \rho \frac{\partial }{{\partial \rho }}\left( {\frac{{1 - {\rho ^{k + 1}}}}{{1 - \rho }}} \right) =\rho \cdot \frac{{1  - {\rho ^k}\left( {1 + k(1 - \rho )} \right)}}{{{{\left( {1 - \rho } \right)}^2}}}.\]
Equation (\ref{P and L relationship-mmk}) is a new relationship. As above, we proved that this relationship holds for \emph{any} subchain 1, not just when subchain 1 is an M/M/s/s queue. We can easily confirm that Equation (\ref{P and L relationship-mmk}) converges to Equation (\ref{P and L relationship-mm1}) at the limit of $k \to \infty $.

3) \textcolor{black}{Let subchain 2 be a queueing system with a harmonic discouragement of arrivals with respect to the number present in the system:} $\lambda _k = \alpha /(k+1), k=0,1,2,...$ and $\mu_k = \mu, k=1,2,3,...$ for subchain 2. Let $\rho=\alpha/\mu$. According to \citet{kleinrock1975queueing}, we know $\pi _s^2 =e^{-\rho}$ and $L_Q^2 = \rho$, from which we obtain
\begin{equation}
\label{discourage sharing a single state}
{P_{block}} = \frac{{e^{-\rho}{P_Q}}}{{1 - \left( {1 - e^{-\rho}} \right){P_Q}}}= \frac{{L_Q}}{{\rho e^{\rho} - \left( {e^{\rho} - 1} \right)L_Q}}.
\end{equation}
Equation (\ref{discourage sharing a single state}) is a new relationship, which holds for \emph{any} subchain 1.

\subsection{Nested Subchains}
Consider a set of \emph{nested} MCs $\{A_k:k=0,1,2,\cdots\}$. $A_k$ is composed of $k+1$ states: ${A_k} = \left\{ {0,1,2,...,k} \right\}$, which satisfies ${A_k} \supset {A_{k - 1}}$ and $\{k\}  = {A_k}\backslash {A_{k - 1}}$ for $\forall k \in {Z^ + }$. Note that $A_0 =\{0\}$. Assume that the boundary condition (Proposition \ref{boundary condition}) is always satisfied throughout the analysis. That is, every time we decompose an MC, we assume that an appropriate termination is applied to each decomposed subchain; hence, the steady-state distribution of the decomposed subchain is always proportional to the original steady-state before the decomposition is made. We are interested in finding a recursive equation for the steady-state probability of state $k$ in ${A_k}$. We decompose $A_k$ into $A_{k-1}$ and $\{k\}$, and apply Theorem \ref{total expectation}. Let state $k$ be the reference state. For any function $f(X)$ of states, the following recursive equation holds:
\begin{equation}
\label{recursive equation}
\frac{{E_{k}\left[ {f\left(X \right)} \right]}}{\pi_k^k} = \frac{E_{k - 1}\left[{f(X)} \right]}{\beta_{k,k-1} \pi_{k-1}^{k - 1}} + f(k), k \in {Z^ + }.
\end{equation}

\remark{By repeating the recursive process and utilizing the properties $\beta_{k,i} \cdot \beta_{i,j}=\beta_{k,j}$ and $\beta_{k,k}=1$, Equation (\ref{recursive equation}) is reduced to 
\begin{equation}
\label{recursive equation2}
\frac{E_k [f(X)]}{\pi_k^k}=\sum_{i=0}^k \frac{f(i)}{\beta_{k,i}}.
\end{equation}
Equation (\ref{recursive equation2}) is immediately obtained by applying Theorem \ref{total expectation} to a set of decomposed ``subchains'' $\{\{0\},\{1\},\cdots,\{k\}\}$, with reference state $k$. Note also that Equation (\ref{recursive equation2}) is equivalent to the definition of expectation: By multiplying $\pi_k^k$ to both sides of the equation, we recover $E_k [f(X)]=\sum_{i=0}^k f(i)\pi_i^k$.}
\textcolor{black}{\remark{One important application of a model with nested queues is a make-to-stock inventory system with priorities. \cite{abouee2012strategies} analyze a \emph{multilevel rationing} (MR) policy for this problem by considering backlog queues, along with a first-come, first-served (FCFS) policy and a \emph{strict priority} (SP) policy (which is a special case of the MR policy). Under the MR policy, customers are segmented to priority levels with different (nondecreasing) threshold inventory levels, where a backlog queue for customers in a given priority level does not impact the backlog queue for higher-priority customers, but only affects the backlog queue for lower-priority customers; such backlog queues with multiple priority classes are effectively represented by a set of nested queues \citep[see \S 3.3 in][]{abouee2012strategies}.}}

Equation (\ref{recursive equation}) is a general property that holds for any finite nested MC and for any function of an MC state. To derive some practically useful properties, we consider a birth and death MC $A_k=\{0,1,2,\cdots,k\}$ with arrival $\lambda_i$ at state $i$ ($i=0,1,\cdots,k-1$) and departure rates $\mu_i$ at state $i$ ($i=1,\cdots,k$). Denote $\rho_{i-1} \doteq \lambda _{i-1} / \mu _{i}$. Denote also that $P_{block}^k \doteq \pi_k^k$ and $L^k \doteq E_k [N]$ for a system $A_k$, where an operator $N$ represents the number of people in the system. (Note that $P_{block}^0 = \pi _0^0 = 1$.) Notice that from the local flow balance equation between states $k-1$ and $k$, we have $\beta _{k,k - 1} = \pi _k / \pi _{k - 1} = \lambda _{k - 1} / \mu _k = \rho _{k - 1}$. By letting $f(X) = 1$ and $N$ in Equation (\ref{recursive equation}), we obtain, correspondingly,
\begin{equation}
\label{recursive equation3}
\frac{1}{P_{block}^k} = \frac{1}{{\rho_{k-1}}P_{block}^{k-1}} + 1 \text{ and } \frac{L^k}{P_{block}^k} = \frac{L^{k - 1}}{{\rho _{k - 1}}P_{block}^{k - 1}} + k.
\end{equation}
Equation (\ref{recursive equation3}) holds for any birth and death MC. To find recursive equations for $P_{block}^k$ and $L^{k}$, we do not need to derive their actual representations; instead, we need concrete representations of $\rho_{k-1}$. In particular,

1) If an MC is an M/M/1/k queue, then by setting $\rho _{k - 1} = \rho  = \lambda/\mu$ for $\forall k \in {Z^ + }$, we obtain \[\frac{1}{{P_{block}^k}} = \frac{1}{{{\rho _{k - 1}}P_{block}^{k - 1}}} + 1 = \frac{1}{{\rho P_{block}^{k - 1}}} + 1 \text{ and } \frac{{{L^k}}}{{P_{block}^k}} = \frac{{{L^{k - 1}}}}{{{\rho _{k - 1}}P_{block}^{k - 1}}} + k = \frac{{{L^{k - 1}}}}{{\rho P_{block}^{k - 1}}} + k.\]
These relationships can be confirmed by plugging in the actual representations of $P_{block}^k$ and $L^k$; to the best of our knowledge, the second relationship is new.

2)  If an MC is an M/M/k/k queue, then by setting ${\rho _{k - 1}} = \lambda / (k\mu )$, we obtain \[\frac{1}{{P_{block}^k}} = \frac{1}{{{\rho _{k - 1}}P_{block}^{k - 1}}} + 1 = \frac{{k\mu }}{{\lambda P_{block}^{k - 1}}} + 1 \text { and }\frac{{{L^k}}}{{P_{block}^k}} = \frac{{{L^{k - 1}}}}{{{\rho _{k - 1}}P_{block}^{k - 1}}} + k = \frac{{k\mu {L^{k - 1}}}}{{\lambda P_{block}^{k - 1}}} + k.\] This first relationship is known as the Erlang B recursive formula (see \citealt{kleinrock1975queueing}). The second is new.

3)  If an MC is an Engset queue, then by setting ${\rho _{k - 1}} = (M - k)\lambda / (k\mu ) $, we obtain  \[\frac{1}{{P_{block}^k}} = \frac{1}{{{\rho _{k - 1}}P_{block}^{k - 1}}} + 1 = \frac{{k\mu }}{{(M-k)\lambda P_{block}^{k - 1}}} + 1 \text { and }\frac{{{L^k}}}{{P_{block}^k}} = \frac{{{L^{k - 1}}}}{{{\rho _{k - 1}}P_{block}^{k - 1}}} + k = \frac{{k\mu {L^{k - 1}}}}{{(M - k)\lambda P_{block}^{k - 1}}} + k.\] The first relationship is known as the Engset recursive formula (see \citealt{kleinrock1975queueing}). The second is new.

\section{Performance Indicators for an M/M/s/k Queue}
\label{sec: EC-mmsk}
We derive properties of interest, $1/\pi_s$ and $L_Q/\pi_s$, when subchains share a single state. A closed-form representation of these quantities is utilized in \S\ref{sec: CBS}. Consider an M/M/s/k queue where $s$ is the number of servers and $k$ is the capacity of the system. Denote the arrival rate as $\lambda$ and the service rate for each server as $\mu$. We decompose the full MC $A=\{0,1,2,\cdots,k\}$ into two subchains sharing a single state $s$: $A_1=\{0,1,\cdots,s\}$ and $A_2=\{s,s+1,\cdots,k-1,k\}$. Since subchains are connected at a single state, truncation is sufficient to conserve their steady-state distributions. Both truncated subchains are well-known queueing systems: $A_1$ is an M/M/s/s queue and $A_2$ is an M/M/1/k-s queue with a utilization parameter $\rho=\lambda/(s\mu)$.

To simplify the representation, let $X$ be a Poisson random variable with parameter $\lambda /\mu$ : $E[X]=var(X)=\lambda /\mu $. Then, we know  (for example, see \citealt{harchol2013performance} and \citealt{kleinrock1975queueing})
$$\frac{1}{\pi_s^1} =\frac{\sum_{i=0}^s \frac{(\lambda/\mu)^i}{i!}}{\frac{(\lambda/\mu)^s}{s!}}=\frac{\Pr \{X \le s\}}{\Pr \{X = s\}} \text{ and } \frac{1}{\pi_s^2}=1+\rho+\cdots+\rho^{k-s}=\frac{1-\rho^{k-s+1}}{1-\rho}.$$
Hence, using Equation (\ref{sharing a single state}), we obtain $1/\pi_k$ and $L_Q/\pi_k$:
$$\frac{1}{\pi_s}=\frac{1}{\pi_s^1}+\frac{1}{\pi_s^2}-1=\frac{\Pr \{X \le s\}}{\Pr \{X = s\}}+\frac{\rho \cdot (1-\rho^{k-s})}{1-\rho}.$$
\begin{eqnarray*}
\frac{L_Q}{\pi_s}&=&\frac{L_Q^2}{\pi_s^2}=\frac{1\cdot \pi_{s+1}^2+2 \cdot \pi_{s+2}^2 + \cdots + (k-s) \cdot \pi_k^2}{\pi_s^2} \\
&=&\frac{1\rho \cdot \pi_s^2+2\rho^2 \cdot \pi_s^2 + \cdots + (k-s)\rho^{k-s} \cdot \pi_s^2}{\pi_s^2}\\
&=&\rho \cdot (1+2\rho+\cdots+(k-s)\rho^{k-s-1})\\
&=&\rho \frac{\partial}{\partial \rho} (1+\rho+\rho^2+\cdots+\rho^{k-s})=\rho\frac{\partial}{\partial \rho}\left(\frac{1-\rho^{k-s+1}}{1-\rho}\right)\\
&=&\frac{\rho \cdot [-(k-s+1)\rho^{k-s}(1-\rho)+(1-\rho^{k-s+1})]}{(1-\rho)^2}\\
&=&\frac{\rho \cdot [1-(k-s+1)\rho^{k-s}+(k-s)\rho^{k-s+1}]}{(1-\rho)^2}.
\end{eqnarray*}

To make use of these results in \S\ref{sec: CBS}, we want to convert the reference state from $s$ to $k$ and the parameter from $\rho$ to  $\omega=1/\rho$. To convert the reference state, notice that $\beta_{k,s}=\pi_k/\pi_s=\rho^{k-s}$. Using this result and parameter $\omega=1/\rho$, we can convert the representation as follows:
$$\frac{1}{\pi_k}=\frac{1}{\beta_{k,s}\pi_s}=\frac{\Pr \{X \le s\}}{\Pr \{X = s\}\rho^{k-s}}+\frac{\rho \cdot (1-\rho^{k-s})}{(1-\rho)\rho^{k-s}}=\frac{\Pr \{X \le s\}\omega^{k-s}}{\Pr \{X = s\}}+\frac{1-\omega^{k-s}}{1-\omega} \doteq f^1(k,\omega).$$
\begin{eqnarray*}
\frac{L_Q}{\pi_k}&=&\frac{L_Q}{\beta_{k,s}\pi_s}=\frac{\rho \cdot [1-(k-s+1)\rho^{k-s}+(k-s)\rho^{k-s+1}]}{(1-\rho)^2 \rho^{k-s}}=\frac{(k-s)-(k-s+1)\omega+\omega^{k-s+1}}{(1-\omega)^2} \doteq g^1(k,\omega).
\end{eqnarray*}

\section{Performance Indicators for an M/M/1 Queue with Restart}
\label{sec: EC-mm1r}
We consider an MC $A_k=\{0,1,2,\cdots,k\}$, which is a birth and death MC with an extra transition (restart) from the last state $k$ back to the first state 0. We are specifically interested in quantities $1/\pi_k^k$ and $L^k/\pi_k^k$. A closed-form representation of these quantities is utilized in \S\ref{sec: CBS}. According to Equation (\ref{recursive equation2}), all we need to know is the expression for $1/\beta_{k,i}(=\pi_i^k/\pi_k^k=\beta_{i,k})$. Let $\lambda_i$ be the arrival rate at state $i$ ($i=0,1,\cdots,k-1$) and $\mu_i$ be the departure rate from state $i$ ($i=1,2,\cdots,k$). Let $r$ be the rate of transition from state $k$ to state 0. Note that the flow balance equation holds: $\lambda_i\pi_i^k=\mu_{i+1}\pi_{i+1}^k+r\pi_k^k$ for $i=0,1,\cdots,k-1$. Denote $a_{i}=\mu_{i+1}/\lambda_i$ and $b_{i}=r/\lambda_i$. Dividing the flow balance equation by $\lambda_i\pi_k^k$, we obtain a recursive equation: $\beta_{i,k}={a_i}\beta_{i+1,k}+b_i$, where $i=0,1,\cdots,k-1$. By applying this recursive equation repeatedly with the condition $\beta_{k,k}=1$, we can obtain the expression for ${\beta_{i,k}}$.

In particular, in \S\ref{sec: CBS}, we consider the simplest case where, for all $i$, $a_{i}=\omega$ (constant) and $b_{i}=1$ hold. To derive $1/\pi_k^k$ and $L^k/\pi_k^k$, we first find $\beta_{i,k}$, which is simply expressed as
$$\beta_{i,k}=1+\omega+\omega^2+\cdots+\omega^{k-i}\beta_{i+(k-i),k}=\frac{1-\omega^{k-i+1}}{1-\omega} \text{, or equivalently, }\beta_{k-i,k}=\frac{1-\omega^{i+1}}{1-\omega}.$$
In addition, we use the following property:
$$\sum_{i=0}^k {i\omega^{i-1}}=\sum_{i=0}^k \frac{\partial \omega^i}{\partial \omega}=\frac{\partial}{\partial \omega} \left(\sum_{i=0}^k \omega^i\right)=\frac{\partial}{\partial \omega} \left(\frac{1-\omega^{k+1}}{1-\omega}\right)=-\frac{(k+1)\omega^k}{1-\omega}+\frac{1-\omega^{k+1}}{(1-\omega)^2}.$$
By letting $f(X)=1$ and $N$ in Equation (\ref{recursive equation2}), we can obtain analytical expressions for $1/\pi_k^k$ and $L^k/\pi_k^k$, respectively:
\begin{eqnarray*}
\frac{1}{\pi_k^k}&=&\sum_{i=0}^k \beta_{i,k}=\sum_{i=0}^k \beta_{k-i,k}=\sum_{i=0}^k \frac{1-\omega^{i+1}}{1-\omega}=\frac{1}{1-\omega}\sum_{i=0}^k 1 - \frac{\omega}{1-\omega}\sum_{i=0}^k \omega^i \\
&=&\frac{k+1}{1-\omega}-\frac{\omega(1-\omega^{k+1})}{(1-\omega)^2}\doteq f^2(k,\omega)
\end{eqnarray*}
and
\begin{eqnarray*}
\frac{L^k}{\pi_k^k}&=&\sum_{i=0}^k i\beta_{i,k} = \sum_{i=0}^k (k-i)\beta_{k-i,k} =\sum_{i=0}^k \frac{(k-i)(1-\omega^{i+1})}{1-\omega}\\
&=& k\sum_{i=0}^k \frac{1-\omega^{i+1}}{1-\omega} - \sum_{i=0}^k \frac{i}{1-\omega} + \frac{\omega^2}{1-\omega} \sum_{i=0}^k i\omega^{i-1}\\
&=& \frac{k(k+1)}{1-\omega}-\frac{k\omega(1-\omega^{k+1})}{(1-\omega)^2} -\frac{k(k+1)}{2(1-\omega)}-\frac{(k+1)\omega^{k+2}}{(1-\omega)^2}+\frac{(1-w^{k+1})\omega^2}{(1-\omega)^3}\\
&=& \frac{k(k+1)}{2(1-\omega)}-\frac{(k+\omega^{k+1})\omega}{(1-\omega)^2}+\frac{(1-\omega^{k+1})\omega^2}{(1-\omega)^3}\\
&=& \frac{k(k+1)}{2(1-\omega)} - \frac{(k-(k+1)\omega+\omega^{k+1})\omega}{(1-\omega)^3} \doteq g^2(k,\omega).
\end{eqnarray*}

\section{Performance Indicators of Subchains in CBS Model}
\label{sec: EC-CBS}
We decompose the full MC representing the CBS model into five (partially overlapping) subchains: $A_1=\{0_{A},1_{A} ,\cdots,n-1_{A}\}$, $A_2=\{n_{A},n+1_{A} ,\cdots,N-1_{A}\}$, $A_3=\{N-1_A, n+1_B\}$, $A_4=\{n+1_{B},n+2_{B} ,\cdots,N_{B}\}$, $A_5=\{N+1_{B},N+2_{B},\cdots\}$. We denote $s=c-e$, $\rho =\lambda/s\mu$, $\omega=1/\rho=s\mu/\lambda$, and $\eta=\lambda/{c\mu}$. We assume that $n, N, c, e$, and $s$ are all integers that satisfy $N>n \geq c$, $e>0$, and $s=c-e>0$. We analyze each subchain independently.

(1) subchain $A_1$: Since this subchain is connected to the rest at a single state, truncation is sufficient to conserve its steady-state distribution. A truncated subchain $A_1$ is a regular M/M/$s$/$k$ queue, whose solution is shown in Appendix \ref{sec: EC-mmsk}.  Let $X$ be a Poisson random variable with parameter $\lambda /\mu :E[X]=var(X)=\lambda /\mu$. By setting $k=n-1$, $s=c-e$, and $\omega=1/\rho=s\mu/\lambda$ for the formulae for $f^1(k,\omega)$ and $g^1(k,\omega)$ in Appendix \ref{sec: EC-mmsk}, we obtain
$$\frac{1}{\pi _{n-1_{A}}^1} =f^1(n-1,\omega) \text{ and } \frac{L_Q^1}{\pi _{n-1_{A} }^1}=g^1(n-1,\omega).$$

(2) subchain $A_2$: Since there is a single inflow state at state $n_A$, we can use Corollary \ref{special termination} to determine the appropriate termination. A terminated subchain $A_2$ is an M/M/1/$k$ queue with restart. Notice that $A_2$ starts from state $n_A$, where $n-s(>0)$ people are already in a queue. Hence, the average number of waiting people in $A_2$ can be obtained by shifting the average number of people in the M/M/1/k queue by $n-s$. By setting $k=(N-1)-n=N-n-1$ and $\omega=1/\rho=s\mu/\lambda$ for the formulae for $f^2(k,\omega)$ and $g^2(k,\omega)$ in Appendix \ref{sec: EC-mm1r}, we obtain

$$\frac{1}{\pi_{N-1_A}^2}=f^2(N-n-1,\omega)$$ and $$\frac{L_Q^2}{\pi_{N-1_A}^2} =\frac{n-s}{\pi_{N-1_A}^2}+g^2(N-n-1,\omega)=(n-s)f^2(N-n-1,\omega)+g^2(N-n-1,\omega).$$

Also, using the expression for $\beta_{i,k}$ in Appendix \ref{sec: EC-mm1r}, by setting $i=0$ and $k=N-n-1$, we obtain
$$\beta_{n_A,N-1_A}=\frac{1-\omega^{N-n}}{1-\omega}.$$

(3) subchain $A_3$: Since there is a single inflow state $n+1_B$ from chain $B$ and a single inflow state $N-1_A$ from chain $A$, we can again use Corollary \ref{special termination} to determine the appropriate termination. A terminated subchain $A_3$ is a two state MC, with a transition from $N-1_A$ to $n+1_B$ at rate $\lambda$ and a transition from $n+1_B$ to $N-1_A$ at a rate $c\mu$. We only need to know the $\beta$ coefficient for this subchain. Using $\eta =\lambda/(c\mu)$,
$$\beta _{N-1_{A},n+1_{B}} =\frac{1}{\eta}.$$

(4) subchain $A_4$: This subchain is symmetric to subchain $A_2$. A terminated subchain $A_4$ is a reversed queueing system of the M/M/1/k queue with restart in Appendix \ref{sec: EC-mm1r}. Hence, the average number of people in $A_4$ can be obtained by subtracting the number in the original M/M/1/k queue with restart from its capacity $k=N-n-1$. As in (2) above, we need to shift the number by $n+1-c(\geq 0)$, who are already in a queue at the left-most state $n+1_B$ in subchain $A_4$. Therefore, we obtain
$$\frac{1}{\pi_{n+1_{B}}^4} =f^2(N-n-1,\eta)$$
and
$$\frac{L_Q^4}{\pi_{n+1_{B}}^4}=\frac{k+(n+1-c)}{\pi_{N-1_A}^2}-g^2(N-n-1,\eta)=(N-c)f^2(N-n-1,\eta)-g^2(N-n-1,\eta).$$

Also, using the expression for $\beta_{i,k}$ in Appendix \ref{sec: EC-mm1r} with the capacity $k=N-n-1$ and the rate $\eta$, we obtain the following expression. (Note that the order of subscript in $\beta$ is reversed because we reverse the numbering of states in the MC in Appendix \ref{sec: EC-mm1r}.)
$$\beta_{N_B,n+1_B}=\frac{1-\eta^{N-n}}{1-\eta}.$$

(5) subchain $A_5$: Since this subchain is connected to the rest at a single state, truncation is sufficient to conserve its steady-state distribution. A truncated subchain $A_5$ is a regular M/M/1 queue with the utilization rate $\eta=\lambda/(c\mu)$, where the solution is well-known. The average waiting people is obtained by shifting the average number by $(N+1)-c$, which is the number of waiting people at the left-most state $N+1_B$ in subchain $A_5$. We obtain:
$$\pi_{N+1_B}^5=1-\eta \text{ and }L_Q^5=N-c+1+\frac{\eta}{1-\eta},$$
or equivalently,
$$\frac{1}{\pi_{N+1_B}^5} =\frac{1}{1-\eta} \text{ and } \frac{L_Q^5}{\pi_{N+1_B}^5} =\frac{N-c+1}{1-\eta} +\frac{\eta}{(1-\eta)^2}.$$
\

The final task is to identify all $\beta$ coefficients. Notice that $\beta_{n_{A},n-1_{A}} =1/\omega$ and $\beta_{N_{B},N+1_{B}} =1/\eta$ hold. Hence, we can derive other necessary coefficients as follows:
$$\beta_{N-1_{A} ,n-1_{A}} =\beta_{N-1_{A} ,n_{A}} \cdot \beta_{n_{A} ,n-1_{A}} =\frac{1-\omega}{\omega (1-\omega^{N-n})} \text{ and}$$
$$\beta_{N-1_{A},N+1_{B}} =\beta_{N-1_{A},n+1_{B} } \cdot \beta_{n+1_{B},N_{B}} \cdot \beta_{N_{B},N+1_{B}} =\frac{1-\eta}{\eta^{2} (1-\eta^{N-n})}.$$

\section{z-Transforms of Subchains in M\textsubscript{t}/M\textsubscript{t}/1 Queue}
\label{sec: EC-MtMt1}
\textcolor{black}{We derive the results indicated in Table \ref{table-MtMt1}. We start our analysis from the $\hat{P}(z)=1$ case and then the $\hat{P}(z)=P(z)$ case.}
\subsection{Outflow from $Q_p$ through a single channel: $\hat{P}(z)=1.$}
\textcolor{black}{We plug $\hat{P}(z)=1$ and $\bar{\gamma}=\sum_{i=0}^\infty \gamma_i p_i =\gamma p_0$ into Equation \eqref{MMPP} and obtain
\begin{equation}
\label{MMPP2}
P(z)=\frac{(1-z)+\frac{{\gamma}}{\mu} z\left[1-\hat{\Pi}(z)\right]}{(1-\frac{\lambda}{\mu}z)(1-z)}p_0.
\end{equation}}
\textcolor{black}{To find $p_0$, we use the condition $P(1)=1$. Since $\hat{\Pi}(1)=1$, $P(1)$ takes the indeterminate form $0/0$. Thus, we must apply L'Hospital's rule, which gives $p_0=\cfrac{1-\frac{\lambda}{\mu}}{1+\frac{\gamma}{\mu}\hat{\Pi}'(1)}$, where $\hat{\Pi}'(z)$ is the first derivative of $\hat{\Pi}(z)$. We now consider two cases for $\hat{\Pi}(z)$.}\\
\textcolor{black}{(i) Inflow into $Q_p$ through a single channel: $\hat{\Pi}(z)=1.$ In this case, both inflow and outflow are observed only at a single state $0_p$ in $Q_p$. Since the termination applied to $Q_p$ is a self-transition $\Delta^p_{0,0}=\gamma_0 \hat{\pi}_{0}=\gamma$ at state $0_p$, our termination scheme becomes a simple truncation, and thus $Q_{\pi}$ does not affect $Q_{p}$. Or we can obtain the same conclusion simply by plugging $\hat{\Pi}(z)=1$ into Equation \eqref{MMPP2}, which is reduced to a familiar form (z-transform of the distribution of M/M/1 queue):
\begin{equation*}
P(z)=\frac{p_0}{1-\frac{\lambda}{\mu}z}, \quad p_0=1-\frac{\lambda}{\mu}.
\end{equation*}}
\textcolor{black}{(ii) Inflow into $Q_p$ through multi channels: $\hat{\Pi}(z)={\Pi}(z).$ This case corresponds to the analysis of $Q_n$ given the knowledge (z-transform of the distribution) of $Q_{n-1}$ in Figure \ref{fig:mtmt1queue}. The termination applied to $Q_p$ is a redirection of outflow from state $0_p$ to all states following $\pi$: $\Delta^p_{0,k}=\gamma_0 \hat{\pi}_{k}=\gamma {\pi}_{k}, \forall k$ (Figure \ref{fig:case-b}), which implies that $Q_{\pi}$ affects $Q_p$. Substituting $\hat{\Pi}(z)$ by 
${\Pi}(z)$ in Equation \eqref{MMPP2} yields
\begin{equation}
\label{MMPP3}
P(z)=\frac{(1-z)+\frac{{\gamma}}{\mu} z\left[1-{\Pi}(z)\right]}{(1-\frac{\lambda}{\mu}z)(1-z)}p_0, \quad p_0=\cfrac{1-\frac{\lambda}{\mu}}{1+\frac{\gamma}{\mu}{\Pi}'(1)},
\end{equation}
where  ${\Pi}'(z)$ is the first derivative of ${\Pi}(z)$ and ${\Pi}'(1)$ is the average queue size of $Q_{\pi}$.}

\subsection{Outflow from $Q_p$ through multi channels: $\hat{P}(z)={P}(z).$}
\textcolor{black}{Plugging $\hat{P}(z)={P}(z)$ and $\bar{\gamma}=\sum_{i=0}^\infty \gamma p_i =\gamma$ into Equation \eqref{MMPP}, we obtain
\begin{equation}
\label{MMPP4}
P(z)=\frac{(1-z)p_0-\frac{\gamma}{\mu}z \hat{\Pi}(z)}{(1-\frac{\lambda}{\mu}z)(1-z)-\frac{\gamma}{\mu}z}.
\end{equation}
We confirm that $P(1)=\hat{\Pi}(1)=1$ satisfies Equation \eqref{MMPP4}. To determine the value of $p_0$, we rewrite Equation \eqref{MMPP4} to explicitly show the poles of $P(z)$ as follows:
\begin{equation*}
P(z)=\frac{(1-z)p_0-\frac{\gamma}{\mu}z \hat{\Pi}(z)}{(1-\frac{z}{r_1})(1-\frac{z}{r_2})},
\end{equation*}
where $r_i, i=1,2$, are the roots of the denominator of $P(z)$. These roots are given by
\begin{equation*}
r_i=\frac{\lambda+\mu+\gamma +(-1)^{i} \sqrt{(\lambda+\mu+\gamma)^2-4\lambda \mu}}{2\lambda},\quad i=1,2.
\end{equation*}
By inspecting the denominator of Equation \eqref{MMPP4} at $z=0,1$, we see that $0<r_1<1,1<r_2$ hold. Since $P(z)$ must be bounded at $|z|<1$, the numerator of $P(z)$ must be zero at $z=r_1$. Thus, we obtain
\begin{equation*}
p_0=\frac{\gamma}{\mu} \frac{ r_1\hat{\Pi}(r_1)}{1-r_1}.
\end{equation*}
We again consider two cases for $\hat{\Pi}(z)$.}\\
\textcolor{black}{(i) Inflow into $Q_p$ through a single channel: $\hat{\Pi}(z)=1$. This case corresponds to the analysis of $Q_1$ (which is independent of $Q_n$) in Figure \ref{fig:mtmt1queue}. Since $\hat{\pi}_0=1$ and $\hat{\pi}_i=0, \forall i\ge1$, the termination applied to $Q_p$ is a simple redirection of all outflows from each state in $Q_p$ to state $0_p$ (including a self-transition at $0_p$): $\Delta^p_{k,0}=\gamma_k \hat{\pi}_0=\gamma, \forall k$ (Figure \ref{fig:case-c}), which implies that $Q_{\pi}$ does not affect $Q_{p}$. (This MC, an M/M/1 queue with extra transitions $\Delta^p_{k,0}=\gamma, \forall k$, is analyzed in \cite{nelson1995probability}, where the model is called \emph{A Processor Model with Failures}.) Substituting $\hat{\Pi}(z)=1$ into Equation \eqref{MMPP4} yields
\begin{equation}
\label{MMPP5}
P(z)=\frac{(1-z)p_0-\frac{\gamma}{\mu}z}{(1-\frac{\lambda}{\mu}z)(1-z)-\frac{\gamma}{\mu}z}, \quad p_0=\frac{\gamma }{\mu} \frac{r_1}{1-r_1},
\vspace*{0.2cm}
\end{equation}
which indeed shows the independence of $P(z)$ from $\pi$ (i.e., $Q_{\pi}$ does not affect $Q_{p}$).}\\
\textcolor{black}{(ii) Inflow into $Q_p$ through multi channel: $\hat{\Pi}(z)={\Pi}(z)$. This final case corresponds to the analysis of $Q_{k}$ for $k=2,3,...,n-1$ in Figure \ref{fig:mtmt1queue}. The termination applied to $Q_p$ is $\Delta^p_{k,k'}=\gamma_k \hat{\pi}_{k'}=\gamma {\pi}_{k'}, \forall k, k'$ (Figure \ref{fig:case-d}), which implies $Q_{\pi}$ affects $Q_{p}$. Plugging $\hat{\Pi}(z)={\Pi}(z)$ in Equation \eqref{MMPP4}, we obtain
\begin{equation}
\label{MMPP6}
P(z)=\frac{(1-z)p_0-\frac{\gamma}{\mu}z {\Pi}(z)}{(1-\frac{\lambda}{\mu}z)(1-z)-\frac{\gamma}{\mu}z}, \quad p_0=\frac{\gamma }{\mu} \frac{r_1}{1-r_1} {\Pi}(r_1).
\vspace*{0.3cm}
\end{equation}}






\end{document}